\documentstyle[12pt,amssymb,amstex]{amsart}

\newtheorem{thm}{Theorem}[section]
\newtheorem{lem}[thm]{Lemma}
\newtheorem{cor}[thm]{Corollary}
\newtheorem{prop}[thm]{Proposition}
\newtheorem{defin}[thm]{Definition}
\newtheorem{rem}[thm]{Remark}

\newtheorem{exam}[thm]{Example}

\newtheorem{conjec}[thm]{Conjecture}
\def\a{\alpha}
\def\e{\varepsilon}
\def\id{{\bf 1}\!\!{\rm I}}

\begin{document}

\title{$G-$Decompositions of Matrices and Related Problems I}

\author{Rasul Ganikhodjaev}
\address{Rasul Ganikhodjaev\\
Faculty of Mechanics and Mathematics, National University of Uzbekistan\\
Vusgorodok, 100174, Tashkent, Uzbekistan}

\author{Farrukh Mukhamedov}
\address{Farrukh Mukhamedov\\
Department of Computational \& Theoretical Sciences Faculty of Sciences\\
International Islamic University  Malaysiya\\
P.O. BOX,141, 25710, Kuantan, Pahang, Malaysiya} \email{{\tt
far75m@@yandex.ru}, {\tt farrukh\_m@@iiu.edu.my}}

\author{Mansoor Saburov}
\address{Mansoor Saburov\\
Department of Computational \& Theoretical Sciences Faculty of Sciences\\
International Islamic University  Malaysiya\\
P.O. BOX,141, 25710, Kuantan, Pahang, Malaysiya} \email{{\tt
msaburov@@gmail.com}}

\begin{abstract}
In the present paper we introduce a notion of $G-$decompositions of
matrices. Main result of the paper is that a symmetric matrix $A_m$
has a $G-$decomposition in the class of stochastic (resp.
substochastic) matrices if and only if $A_m$ belongs to the set
${\mathbf{U}}^m$ (resp. ${\mathbf{U}}_m$). To prove the main result,
we study extremal points and geometrical structures of the sets
${\mathbf{U}}^m$, ${\mathbf{U}}_m$. Note that such kind of
investigations enables to study Birkhoff's problem for quadratic
$G-$doubly stochastic operators.
\end{abstract}

\maketitle

\section{Introduction}

Let us recall that a matrix $A_m=(a_{ij})_{i,j=1}^m$ is said to be
\begin{enumerate}
  \item [(i)] {\it stochastic} if its elements are non-negative and each row sum is equal to one;
  \item [(ii)] {\it substochastic} if its elements are non-negative and each row sum is less or equal to one;
  \item [(iii)] {\it doubly stochastic} if its elements are non-negative and each row and column sums are equal to one.
\end{enumerate}
In \cite{B} G.D Birkhoff characterized the
set of extreme doubly stochastic matrices. Namely his result states
as follows: the set of extreme points of the set of $m\times m$
doubly stochastic matrices coincides with the set of all
permutations matrices.

One can consider a generalization of Birkhoff's result in two
directions. In the first direction, one may consider the description
of all extreme points of the set of \emph{infinite doubly stochastic
matrices}, and in the second one, one may consider the description of
all extreme points of the set of \emph{nonlinear doubly stochastic
operators}.

Concerning the first case, in \cite{K}, \cite{Mld}, the Birkhoff's
problem have been solved, i.e. it was proved that there are no
extreme points of the set of all infinite doubly stochastic matrices
except the permutation matrices. In \cite{S2005,S2006} Yu. Savarov
has shown that, under certain conditions, Birkhoff's result on
doubly stochastic matrices remains valid for countable families of
discrete probability spaces which have nonempty intersections. Let
us also mention some other  related results. For example, in
\cite{Mir} it was proved that an extreme doubly substochastic matrix
is a subpermutation matrix. For its generalization to arbitrary
marginal vectors see \cite{Bru}, for the finite dimensional case and
\cite{Gr2,P}, for the infinite dimensional case. In
\cite{Katz1},\cite{Katz2} the extreme symmetric stochastic and
substochastic matrices, respectively, were determined. These results
were generalized to finite symmetric matrices with given row sums by
R. A. Brualdi \cite{Bru}. Finally in \cite{Ckatz},\cite{Gr1} the
extreme points of the set of infinite symmetric stochastic matrices
with given row sums were described.

The present paper is related to the Birkhoff's problem for nonlinear
doubly stochastic operators. In the this case, we will face with a
few contretemps. In fact, first of all, we should define a
conception of stochasticity for nonlinear operators. We then should
define doubly stochasticity of nonlinear operators. After all of
these, we can consider Birkhoff's problem for nonlinear operators.
However, a conception of doubly stochasticity for nonlinear
operators can be given by different ways. Here, we shall present one
of conceptions of doubly stochasticity in nonlinear settings
introduced in \cite{G}.

Let us recall some necessary notions and notations.

Let $I_m=\{1,2,\cdots,m\}$ be a finite set and $S^{m-1}$ be an
$m-1$ dimensional simplex, i.e.,
$$S^{m-1}=\left\{x=(x_1,x_2,\cdots,x_m)\in {\mathbb{R}}^m: \ \sum\limits_{i=1}^mx_i=1,\quad x_i\geq0\right\}.$$
Every element of the simplex $S^{m-1}$ can be considered as a
probability distribution of the finite set $I_m.$ Hence, the simplex
$S^{m-1}$ is a set of all probability distributions of the finite
set $I_m.$

Any operator $V$ which maps the simplex $S^{m-1}$ into
itself is called a \emph{stochastic} operator.

For a vector $x=(x_1,x_2,\cdots,x_m)\in{\mathbb{R}}^m$ we denote by
$x_{\downarrow}=(x_{[1]},x_{[2]},\cdots,x_{[m]})$ the vector with
same coordinates, but sorted in non-increasing order $x_{[1]}\geq
x_{[2]}\geq \cdots\geq x_{[m]}.$ For $x,y\in{\mathbb{R}}^m,$ we say
that $y$ is majorized by $x$ (or $x$ majorizes $y$), and write
$y\prec x$ if
\begin{eqnarray*}
\sum\limits_{i=1}^ky_{[i]}\leq\sum\limits_{i=1}^kx_{[i]}, \
\textrm{for  all} \ k=\overline{1,m-1},\ \textrm{ and} \
\sum\limits_{i=1}^my_{[i]}=\sum\limits_{i=1}^mx_{[i]}.
\end{eqnarray*}
The Hardy-Littlewood-Polya theorem (see \cite{MO}) says that $y$ is
majorized by $x$, i.e., $y\prec x$ if and only if there exists a
doubly stochastic matrix $A_m$ such that $y=A_mx.$ As a corollary we
can get that a matrix $A_m$ is a doubly stochastic if and only if
$A_mx\prec x$ for any $x\in {\mathbb{R}}^m.$ Thus, we can give
another equivalent definition of the doubly stochasticity of the
matrix as follows: a matrix $A_m$ is called {\it doubly stochastic}
if $A_mx\prec x$ for any $x\in {\mathbb{R}}^m.$

Based on this result, in \cite{G} it has been introduced a
definition of doubly stochasticity for nonlinear operators. Namely,
a stochastic operator $V:S^{m-1}\to S^{m-1}$ is called
$G-$\emph{doubly stochastic} if $Vx\prec x$ for any $x\in S^{m-1}.$
An advantage of this definition is that $G-$doubly stochastic
operators are well defined for any kind of nonlinear stochastic
operators, even though the forms of nonlinear operators are not
polynomial There is another way to define the notion of doubly
stochasticity for quadratic stochastic operators.

Among nonlinear operators, the simplest one is a quadratic one. Such
a quadratic operator $V:{\mathbb{R}}^{m}\to {\mathbb{R}}^{m}$ can be
given as follows
$$Vx=\left(\sum\limits_{i,j=1}^mA_{ij,1}x_ix_j,\sum\limits_{i,j=1}^mA_{ij,2}x_ix_j,\cdots, \sum\limits_{i,j=1}^mA_{ij,m}x_ix_j\right),$$
where ${\mathbb{A}}_V=(A_{ij,k})_{i,j,k=1}^m$ is a cubic matrix. One
can see that every quadratic operator is uniquely defined by a cubic
matrix ${\mathbb{A}}_V$. In fact, if we denote
$A^{(k)}_m=(A_{ij,k})_{i,j=1}^m$ then the quadratic operator has the
following form
$$Vx=\left((A^{(1)}_mx,x),\cdots, (A^{(m)}_mx,x)\right),$$
where, $(\cdot,\cdot)$ is the standard inner product in
${\mathbb{R}}^n$.

In what follows, we shall use the notation
$\left(A^{(1)}_m\mid\cdots\mid A^{(m)}_m\right)$ for the quadratic operator $V$.

 In this paper we attempt to deal with
 Birkhoff's problem for quadratic $G-$doubly stochastic operators \footnote{Here for
the sake of completeness we should mention that there is also
another way to define quadratic doubly stochastic operators is the
following sense: a quadratic operator $V$ is called $Z-$\emph{doubly
stochastic} if its cubic matrix ${\mathbb{A}}_V$ satisfies the
following conditions
$\sum\limits_{i=1}^mA_{ij,k}=\sum\limits_{j=1}^mA_{ij,k}=\sum\limits_{k=1}^mA_{ij,k}=1,$
$A_{ij,k}\geq 0$, for all $ i,j,k=\overline{1,m}.$ One can easily
check that if $V$ is a quadratic $Z-$doubly stochastic then $V$ is
stochastic. Note that $Z-$doubly stochasticity of quadratic
operators differs from $G-$doubly stochasticity. However, the
disadvantage of the this definition is that $Z-$doubly stochastic
operators are only well defined for polynomial nonlinear stochastic
operators. In \cite{Li} it was concerned with possible
generalizations of Birkhoff's problem to higher dimensional
stochastic matrices and provided lots of criteria for extremity of
such matrices. However, the provided criteria given in \cite{Li} is
difficult to check in practice. Therefore, up to now, there is no a
full explicitly description of extreme higher dimensional stochastic
matrices. Particulary, there is not a full explicit description of
extreme quadratic $Z-$doubly stochastic operators as well. In
\cite{SMSF}, it was checked one class of quadratic $Z-$doubly
stochastic operators to be extreme.}.

Let us define the following sets
\begin{eqnarray*}
&&\textbf{U}_m = \left\{A_m=(a_{ij})_{i,j=1}^m: a_{ij}=a_{ji}\ge 0 ,
\ \sum\limits_{i,j\in \a}a_{ij}\le |\a|, \ \forall \a\subset
I_m\right\},\\
&&\textbf{U}^m = \left\{A_m\in \textbf{U}_m: \
\sum\limits_{i,j=1}^ma_{ij}=m\right\},
\end{eqnarray*}
where $|\a|$ stands for a number of elements of a set $\a.$

In \cite{GS}, the investigation of extreme quadratic $G-$doubly
stochastic operators has been started. One of the main results of
the paper \cite{GS} is that if a quadratic stochastic operator
$V=\left(A^{(1)}_m\mid\cdots\mid A^{(m)}_m\right)$ is $G-$doubly
stochastic, then the corresponding $m\times m$ matrices $A^{(k)}_m$
belong to the set $\textbf{U}^m$ for any $k=\overline{1,m}.$ In
other words, the set of all quadratic $G-$doubly stochastic
operators is a convex subset of the set
$\underbrace{\textbf{U}^m\times\cdots\times\textbf{U}^m}_{m}.$ It is
clear that the sets $\textbf{U}^m$ and
$\underbrace{\textbf{U}^m\times\cdots\times\textbf{U}^m}_{m}$ are
convex. Before studying extreme points of the set of all quadratic
$G-$doubly stochastic operators, it is of independent interest to
study geometrical structures of $\textbf{U}^m$. The relationship
between extreme points of the set of all quadratic $G-$doubly
stochastic operators and extreme points of $\textbf{U}^m$ is given
in \cite{GS}: let $V=\left(A^{(1)}_m\mid\cdots\mid A^{(m)}_m\right)$
be a quadratic $G-$doubly stochastic operator. If any $m-1$ matrices
of the matrices $\{A_m^{(k)}\}_{k=1}^m$ are extreme in the set
$\textbf{U}^m$ then the corresponding quadratic $G-$doubly
stochastic operator $V=\left(A^{(1)}_m\mid\cdots\mid
A^{(m)}_m\right)$ is extreme in the set of all quadratic $G-$doubly
stochastic operators. This result encourages us to study extreme
points of $\textbf{U}^m$.

One of the crucial point in the Birkhoff's problem for quadratic
$G-$doubly stochastic operators is a notion of $G-$decomposition of
symmetric matrices. Namely, let ${\mathcal{M}}_{m\times m}$ be a set
of all $m\times m$ matrices, and ${\mathcal{G}}\subset
{\mathcal{M}}_{m\times m}$ be a convex bounded polyhedron.

\begin{defin}
We say that a matrix $A_m$ has a $G-$decomposition in a class
$\mathcal{G}$ if there exists a matrix $X_m\in \mathcal{G}$ such
that
\begin{eqnarray}\label{G-decom}
A_m=\frac{X_m+X^{t}_m}{2},
\end{eqnarray}
By ${\mathcal{G}}^{s}$ we denote the class of all such kind of
matrices $A_m.$ The set ${\mathcal{G}}^{s}$ is called the
symmetrization of $\mathcal{G}.$
\end{defin}

Note that such a notion of $G-$decomposition is related to certain
problems in the control theory.\footnote{Recall that the Lyapunov
equation has a form $Y_mX_m+X_m^{t}Y_m=A_m$ which appears in many
branches of the control theory, such as stability analysis and
optimal control. Here $A_m,X_m$ are $m\times m$ given matrices and
$Y_m$ is an $m\times m$ unknown matrix. If one considers the
Lyapunov equation with respect to $X_m$ and $Y_m=\frac12 \id_{m}$,
where $\id_m$ is the unit matrix, then we get \eqref{G-decom}. Some
observations show that if a symmetric matrix $A_m$ has the
decomposition \eqref{G-decom} for the special matrix $X_m$ then
certain problems of convex analysis can be easily solved.} One of
the fascinating result of the paper \cite{GS} is the following one.

\begin{thm}\label{GanikhodjaevTheorem}
\cite{GS}
Let $A_m$ be a symmetric matrix. Then the following statements are equivalent:
\begin{enumerate}
  \item [(i)] The matrix $A_m$ belongs to the set ${\mathbf{U}}^m$;
  \item [(ii)] The matrix $A_m$ has a $G-$decomposition in class of stochastic matrices;
  \item [(iii)] The inequality $x_{[m]}\le(A_mx,x)\le x_{[1]}$ holds for all $x\in S^{m-1}$.
\end{enumerate}
\end{thm}

To be fair, we would say that in the paper \cite{GS} the provided
proof of the part $(i)\Leftrightarrow (ii)$ of Theorem
\ref{GanikhodjaevTheorem} had some gaps. To clarify and fill those
gaps, we aim to write the present paper as a complementary one to
\cite{GS}. Here, we are going to give a complete proof of Theorem
\ref{GanikhodjaevTheorem}. Moreover, we shall generalize it for
substochastic matrices as well. As we already mentioned above there
is a relationship between extreme points of the set $\textbf{U}^m$
and the set of quadratic $G-$doubly stochastic operators. However,
the extreme points of the set $\textbf{U}^m$ were not described in
\cite{GS}. Therefore, we are going to deeply study geometrical and
algebraical structures of the sets $\textbf{U}^m$ and
$\textbf{U}_m$. This paper contains many results which are of
independent interest.

Let us briefly explain the organization of the paper. The main
results of this paper is the following theorem.

\begin{thm}\label{mainresults} The following statements hold true:
\begin{enumerate}
             \item [(i)] A symmetric matrix $A_m$ has a $G-$decomposition in the class of
stochastic matrices if and only if $A_m$ belongs to the set ${\mathbf{U}}^m$;
             \item [(ii)] A symmetric matrix $A_m$ has a $G-$decomposition in the class of
substochastic matrices if and only if $A_m$ belongs to the set ${\mathbf{U}}_m$.
\end{enumerate}
\end{thm}

The strategy of the proof of the main results is the following: in
both cases it is enough to prove the assertions of the theorem for
extreme points of the sets ${\mathbf{U}}_m$ and ${\mathbf{U}}^m$.
Then we shall employ the Krein-Milman's theorem to prove the theorem
in general setting. First of all, we shall prove the case $(i)$ for
extreme points of ${\mathbf{U}}^m$, then using canonical forms of
extreme points of ${\mathbf{U}}_m$ we reduce the case $(ii)$ to
$(i)$.  Therefore, our first task is to study extreme points of the
sets ${\mathbf{U}}_m$ and ${\mathbf{U}}^m$.

In section \ref{section2} we shall provide criteria for extreme
points of the sets ${\mathbf{U}}_m$ and ${\mathbf{U}}^m$. We stress
here, that the implication $(i)\Rightarrow (iii)$ of Theorem
\ref{localcriterionforU^m} was stated in the paper \cite{GS} without
any justification. Actually, this implication was a main point of
the implication $(i)\Leftrightarrow (ii)$ of Theorem
\ref{GanikhodjaevTheorem}. In this section, we shall justify it,
moreover we show that the inverse implication $(iii)\Rightarrow (i)$
of Theorem \ref{localcriterionforU^m} is also valid. The main
results of this section are Theorems \ref{localcriterion} and
\ref{localcriterionforU^m}.

In section \ref{section3} we shall study explicit forms of the
extreme points of the sets ${\mathbf{U}}_m$ and ${\mathbf{U}}^m$,
respectively. The results of this section would be used to solve
Birkhoff's problem for quadratic $G-$doubly stochastic operators.
The main results of this section are Corollaries \ref{bestcriterion}
and \ref{bestcriterionforU^m}.

In section \ref{section4} we shall study canonical forms of the
extreme points of ${\mathbf{U}}_m$  which is an extremely important
to prove the case $(ii)$ of Theorem \ref{mainresults}. Using the
canonical forms of the extreme points of ${\mathbf{U}}_m$ we are
able to reduce the case $(ii)$ of Theorem \ref{mainresults} to the
case $(i)$ of the same theorem. The main results of this section is
Corollary \ref{inclusionsExtrUm}

In section \ref{section5} we shall prove the main results of this
paper. They are provided by Theorems \ref{stochasticcase} and
\ref{substochasticcase}. There, by means of the results of section
\ref{section2} we first prove Theorem \ref{stochasticcase} then
again using the results of section \ref{section4} we reduce the
proof of Theorem \ref{substochasticcase} to Theorem
\ref{stochasticcase}.

\section{Some criteria for extreme points of the sets $\textbf{U}_m$ and $\textbf{U}^m$}\label{section2}

In this section we want to give some criteria for extreme points of
the sets $\textbf{U}_m$ and $\textbf{U}^m$. Moreover, we provide a
proof of some facts which were not proven in \cite{GS}.

It is clear that $\textbf{U}_m$ is a convex set and $\textbf{U}^m$
is a convex subset of $\textbf{U}_m$.

One can easily see that
\begin{eqnarray*}
\textbf{U}_1 \hookrightarrow \textbf{U}_2 \hookrightarrow \textbf{U}_3 \hookrightarrow \cdots \hookrightarrow \textbf{U}_m.
\end{eqnarray*}
Here, the inclusion $\textbf{U}_k \hookrightarrow \textbf{U}_m $
should be understood in the way that the matrix with smaller order
can be extended to larger by letting new entries to be zero. More
precisely, the inclusion $\textbf{U}_k \hookrightarrow \textbf{U}_m
$ means that if $A_k\in \textbf{U}_k$ then there exists $A_m\in
\textbf{U}_m$ such that
\begin{eqnarray*}
A_m=\left(
      \begin{array}{cc}
        A_k & \circleddash_{k\times m-k} \\
        \circleddash_{m-k\times k} & \circleddash_{m-k\times m-k} \\
      \end{array}
    \right),
\end{eqnarray*}
where $\circleddash_{m\times n}$ means a $m\times n$ matrix with
zero entries.

In the same way, we can get that
\begin{eqnarray*}
\textbf{Extr}\textbf{U}_1 \hookrightarrow \textbf{Extr}\textbf{U}_2 \hookrightarrow \textbf{Extr}\textbf{U}_3 \hookrightarrow \cdots \hookrightarrow \textbf{Extr}\textbf{U}_m,
\end{eqnarray*}
here, $\textbf{Extr}\textbf{U}_k$ denotes the set of the
extreme points of $\textbf{U}_k.$

We are going to study a geometrical structure of the
set $\textbf{U}_m$. Particularly, we describe extreme
points of $\textbf{U}_m$. Let us recall some well-known notations.

A submatrix $A_{\a}$ of $A_m=(a_{ij})_{i,j=1}^m$ is said to be
\emph{a principal submatrix} if $A_{\a}=(a_{ij})_{i,j\in\a}$, i.e.
all entries indexes of $A_{\a}$ belong to $\a(\subset I_m).$

The proof of the following proposition is straightforward.

\begin{prop} The following statements hold true:
\begin{itemize}
  \item [(i)] If $A_m\in {{\mathbf{U}}_m},$ then its any principal submatrix of
order $k$ belong to ${\mathbf{U}}_k,$ $1\le k\le m$;
  \item [(ii)]\label{extremeprinsubmatix} If a principal submatrix $A_{\a}$ of $A_m$ is extreme in
${\mathbf{U}}_{|\a|}$, for some $\a\subset I_m,$ then for any matrices
$A_m', A_m''\in {\mathbf{U}}_m$ satisfying $2A_m=A_m'+A_m'',$ one has
$A_{\a}=A_{\a}'=A_{\a}''.$
\end{itemize}
\end{prop}

\begin{prop}\label{generalcriteraonAm}
A matrix $A_m$ belongs to ${\mathbf{Extr}}{\mathbf{U}}_m$ if and only if for any
entry $a_{ij}$ of $A_m$, there exists $\a\in I_m$ such that $a_{ij}$
is an entry of the principal submatrix $A_{\a}$ with $A_{\a}\in
\mathbf{Extr}\mathbf{U}_{|\a|}.$
\end{prop}

\begin{pf}
\emph{If part.} Let us assume that $A_m$ is not extreme, that is
$2A_m=A_m'+A_m''$, for some $A_m', A_m''\in \textbf{U}_m$ with
$A_m\neq A_m',$ $A_m\neq A_m''.$ From the former, we conclude that
there is an entry $a_{i_0j_0}$ such that $a_{i_0j_0}\neq
a_{i_0j_0}'.$ According to the condition, there exists an extremal
principal submatrix $A_{\a_0}$ in $\textbf{U}_{|\a_0|},$ containing
$a_{i_0j_0}.$ Then due to Proposition \ref{extremeprinsubmatix} we get
$A_{\a_{0}}=A_{\a_0}'=A_{\a_0}'',$ and hence
$a_{i_0j_0}=a_{i_0j_0}'$ which is a contradiction.

\emph{Only If part.} Suppose that $A_m\in
\textbf{Extr}\textbf{U}_m.$ By putting $\a=I_m$ we get $a_{ij}\in
A_{\a}$ and $A_{\a}\in \textbf{Extr}\textbf{U}_{|\a|},$ for any
entry $a_{ij}.$
\end{pf}

\begin{rem} Note that the provided criterion is somehow
difficult to apply in practice. The reason is that sometimes a given
matrix may not have extreme proper principal submatrices. Let us
consider the following $m\times m$ matrix
\begin{equation*}
N_m=\left(
      \begin{array}{cccccccc}
        0 & \frac12 & 0 & 0 & \cdots & 0 & 0 & \frac12 \\
        \frac12 & 0 & \frac12 & 0 & \cdots & 0  & 0 & 0 \\
        0 & \frac12 & 0 & \frac12 & \cdots & 0 & 0 & 0 \\
        \cdots & \cdots & \cdots & \cdots & \cdots & \cdots & \cdots & \cdots \\
        0 & 0 & 0 & 0 & \cdots & \frac12 & 0 & \frac12 \\
        \frac12 & 0 & 0 & 0 & \cdots & 0 & \frac12 & 0 \\
      \end{array}
    \right).
\end{equation*}
For this matrix, the problem of finding of its extreme proper
principal submatrix coincides with the problem of showing its
extremity. Further, one can show that the matrix $N_m$ is not
extreme.
\end{rem}

However, there are some benefits of
the provided criterion in terms of studying some properties of extreme matrices.
The following corollaries directly follow from Proposition \ref{generalcriteraonAm}.

\begin{cor}\label{nonextremeonAm}
A matrix $A_m$  is not extreme in ${\mathbf{U}}_m$ if and only if
there exists an entry $a_{i_0j_0}$ such that any principal submatrix
$A_{\a}$ containing this entry, is not extreme in
$\mathbf{U}_{|\a|}.$
\end{cor}

\begin{cor}\label{2by2matrixextreme} If every $2\times 2$
principal submatrix of a matrix $A_m$ is extreme in ${\mathbf{U}}_2$
then the matrix $A_m$  itself is extreme in ${\mathbf{U}}_m$.
\end{cor}

\begin{rem} The converse of
Corollary \ref{2by2matrixextreme} is not true. For instance, the matrix
\begin{equation*}
M_3=\left(
  \begin{array}{ccc}
    0 & \frac12 & \frac12 \\
    \frac12 & 0 & 0 \\
    \frac12 & 0 & 1 \\
  \end{array}
\right)
\end{equation*}
is extreme in ${\mathbf{U}}_3,$ however it has a $2\times 2$ non
extreme principal submatrix $\left(
  \begin{array}{cc}
    0 & \frac12  \\
    \frac12 & 0  \\
  \end{array}
\right)$.
\end{rem}

Let us now present some facts. We are not going to prove them because of their evidence.

\begin{prop}\label{trivialfacts} The following assertions hold true:
\begin{itemize}
  \item [{\rm(i)}] \label{0aij1} If $A_m\in {\mathbf{U}}_m,$ then $0\le a_{ij}\le 1$ for all $i,j=\overline{1,m};$
  \item [{\rm(ii)}] \label{aij0or1} Let matrices $A_m, A_m', A_m''\in {\mathbf{U}}_m$ satisfy $2A_m=A_m'+A_m''.$
If $a_{i_0j_0}=1\vee 0$ for some $i_0,j_0\in I_m,$ then $a_{i_0j_0}=a_{i_0j_0}'=a_{i_0j_0}''.$
Here and henceforth $a=b\vee c$ means either $a=b$ or $a=c;$
  \item [{\rm(iii)}] \label{either1or0} Any matrix $A_m\in {\mathbf{U}}_m$ having entries being equal to either 1 or 0 is extreme in ${\mathbf{U}}_m$.
\end{itemize}
\end{prop}

Let us introduce the following useful conception.

\begin{defin}
Let $A_m\in {\mathbf{U}}_m$. An index set $\a\subset I_m$ is said to
be saturated, for the matrix $A_m$, whenever $\sum\limits_{i,j\in
\a}a_{ij}=|\a|.$ A principal submatrix $A_\a,$ corresponding to the
saturated index set $\a,$ is called a saturated principal submatrix.
\end{defin}

An advantage of the given conception is that using induction with
respect to the number of saturated index sets we can easily prove
lots of properties of ${\textbf{U}}_m$. Moreover, it is an
appropriate conception to formulate some facts regarding extreme
matrices of ${\textbf{U}}_m.$

\begin{prop}\label{unionandinterofsatutset}
Let $\a,\beta\subset I_m$ be saturated index sets of $A_m\in
{\mathbf{U}}_m.$ Then the following assertions hold true:
\begin{itemize}
  \item[\rm(i)] If $\a\cap\beta\neq\emptyset,$ then
  $\a\cap\beta$ is a saturated index set for $A_m;$
  \item[\rm(ii)] $\a\cup\beta$ is a saturated index set for $A_m$.
\end{itemize}
\end{prop}
\begin{pf}
Let $\gamma=\a\cap\beta$, for any index sets $\a, \beta\subset I_m.$
Then, one can easily check the following equality
\begin{eqnarray}\label{sumalphacupbetaaij}
\sum\limits_{i,j\in\a\cup\beta}a_{ij} =\sum\limits_{i,j\in\a}a_{ij}+
\sum\limits_{i\in\a\setminus\gamma\atop
j\in\beta\setminus\gamma}a_{ij}
+\sum\limits_{i\in\beta\setminus\gamma\atop
j\in\a\setminus\gamma}a_{ij}+\sum\limits_{i,j\in\beta}a_{ij}-\sum\limits_{i,j\in\gamma}a_{ij}.
\end{eqnarray}

Now suppose that $\a,\beta\subset I_m$ are saturated index sets of a
matrix $A_m\in {\mathbf{U}}_m$ and $\gamma\neq\emptyset.$ We only
consider a case when $\gamma\neq\a$, $\gamma\neq\beta$ otherwise the
theorem is evident. It is clear that
$\sum\limits_{i,j\in\gamma}a_{ij}\leq|\gamma|.$ Hence by means of
\eqref{sumalphacupbetaaij} we have
\begin{eqnarray*}\label{sumalphacupbetaaijinequl}
|\a\cup\beta|=|\a|+|\beta|-|\gamma|&\leq&\sum\limits_{i,j\in\a}a_{ij}+\sum\limits_{i,j\in\beta}a_{ij}-\sum\limits_{i,j\in\gamma}a_{ij}\\[2mm]
&\le&
\sum\limits_{i,j\in\a\cup\beta}a_{ij}\le|\a\cup\beta|=|\a|+|\beta|-|\gamma|.
\end{eqnarray*}
This yields that
\begin{eqnarray*}
\sum\limits_{i,j\in\a\cup\beta}a_{ij}=|\a\cup\beta|, \quad
\sum\limits_{i,j\in\gamma}a_{ij}=|\gamma|.
\end{eqnarray*}
 Therefore $\a\cap\beta$
and $\a\cup\beta$ are saturated index sets for $A_m\in
{\mathbf{U}}_m$.
\end{pf}

According to Proposition \ref{unionandinterofsatutset}, a class of
all the saturated index sets of a given matrix is closed with
respect to the operations of union and intersection. That is why,
this property of the saturated index sets implies a reason to
introduce the following

\begin{defin} A saturated principal submatrix of
a matrix $A_m\in {\mathbf{U}}_m$ containing an entry $a_{ij}$ is
called saturated neighborhood of $a_{ij},$ and  the order of such a
saturated principal submatrix is said to be its radius. Saturated
neighborhoods of $a_{ij}$ with minimal and maximal radiuses are
called a minimal and a maximal saturated neighborhoods of $a_{ij},$
respectively.
\end{defin}

\begin{rem} If an entry $a_{ij}$ of $A_m\in {\mathbf{U}}_m$ has a minimal or a maximal
saturated neighborhoods, then they are uniquely defined. Let us show
uniqueness of the minimal saturated neighborhood of $a_{ij}$. Assume
that there are two minimal saturated neighborhoods $A_{\a},$
$A_{\a'}$ of $a_{ij}$, and we denote the corresponding saturated
index sets by $\a$ and $\a'$. Since $\a\cap\a'\neq\emptyset$ and
$A_{\a},$ $A_{\a'}$ are minimal saturated neighborhoods, due to
Proposition \ref{unionandinterofsatutset} (i), $\a\cap\a'$ is a
saturated index set, and  $|\a|\le|\a\cap\a'|,$
$|\a'|\le|\a\cap\a'|.$ Therefore, $\a=\a\cap\a'=\a'.$

Using the same argument with Proposition
\ref{unionandinterofsatutset} (ii), one can get the uniqueness of
the maximal saturated neighborhood of $a_{ij}$.
\end{rem}

We would like to emphasize that the minimal saturated neighborhood plays
an important role, for geometrical structures of the set
${\mathbf{U}}_m$ whereas the maximal
saturated neighborhood plays as crucial point for its algebraical structures.

If an entry $a_{ij}$ has a saturated neighborhood then, since the
matrix $A_m\in {\mathbf{U}}_m$ is symmetric,  an entry  $a_{ji}$ has
also the same saturated neighborhood. That is why, henceforth, we
only consider saturated neighborhoods of entries $a_{ij}$ in which
$i\le j.$

Let us observe the following: assume that $A_\a$ is a principal
submatrix of a matrix $A_m\in {\mathbf{U}}_m$ and $a_{ij}$ is a
entry of $A_\a.$ Then, in general, the minimal saturated
neighborhood of $a_{ij}$ in the matrix $A_m$ does not coincide with
its minimal saturated neighborhood  in the principal submatrix
$A_{\a}.$ For this, one of the main reasons is that the entry
$a_{ij}$ may have a minimal saturated neighborhood in $A_m$, but may
not so in $A_\a.$ We can see this picture in the following example:
let
\begin{equation*}
A_6=\left(
      \begin{array}{cccccc}
        0 & \frac12 & 0 & 0 & 0 & 0  \\
        \frac12 & 0 & \mathbf{\frac12} & 0 & 0  & 0  \\
        0 & \mathbf{\frac12} & 0 & \frac12 & 0 & 0  \\
        0 & 0 & \frac12 & 1 & 0 & 0  \\
        0 & 0 & 0 & 0 & 0 &  0 \\
        0 & 0 & 0 & 0 & 0  & 1\\
      \end{array}
    \right).
\end{equation*}
If we consider the principal submatrix
$A_{\a_0}=\left(\begin{array}{ccc}
0 & \frac12 & 0 \\
\frac12 & 0 & \mathbf{\frac12} \\
0 & \mathbf{\frac12} & 0 \\
\end{array}\right),$
where $\a_0=\{1,2,3\},$ then the element $a_{23}=\frac12$ does not
have a minimal saturated neighborhood in $A_{\a_0}$, but it has the
minimal saturated neighborhood $A_{\beta_0}=\left(\begin{array}{ccc}
0 & \mathbf{\frac12} & 0 \\
\mathbf{\frac12} & 0 & \frac12 \\
0 & \frac12 & 1 \\
\end{array}\right)$
in the given matrix $A_6,$ where $\beta_0=\{2,3,4\}.$

Fortunately, if a principal submatrix $A_\a$ of a matrix $A_m$ is
saturated, then a minimal saturated neighborhood of any element
$a_{ij},$ where $i,j\in\a,$ in the matrix $A_m$, coincides with its
minimal saturated neighborhood in the principal submatrix $A_{\a}.$
Namely, we have the following

\begin{lem}\label{indepenofminimsaturneigh}
Let $A_{\a}$ be a saturated principal submatrix of a matrix
$A_m\in{\mathbf{U}}_m$ and $a_{ij}$ be any entry of $A_\a.$ Then
there exists a minimal saturated neighborhood of $a_{ij}$ in  $A_m$,
and it is a principal submatrix of $A_{\a}.$ Moreover, it is a
minimal saturated neighborhood of $a_{ij}$ in $A_\a.$
\end{lem}

\begin{pf} Suppose that $A_{\a}$ is a saturated principal submatrix of
$A_m\in{\mathbf{U}}_m$ and $a_{ij}$ is an entry of $A_\a.$ Then
$a_{ij}$ has at least one saturated neighborhood in $A_m$, which is
 $A_\a.$ Therefore, it has a minimal saturated neighborhood in
$A_m.$ By $\a'$ we denote a saturated index set corresponding to its
minimal saturated neighborhood in $A_m.$ We must show that
$\a'\subset\a.$ Indeed, since $\a\cap\a'\supset\{i,j\}$ and both
$\a$ and $\a'$ are saturated index sets, then according to
Proposition \ref{unionandinterofsatutset} the set $\a\cap\a'$ is
saturated as well. Further, since $\a'$ corresponds to the minimal
saturated neighborhood of $a_{ij}$ in $A_m$, it follows that
$\a'\subset\a'\cap\a\subset\a$ as desired. Moreover, it immediately
follows from the definition of the minimal saturated neighborhood
and $\a'\subset\a$ that the minimal saturated neighborhood of
$a_{ij}$, in $A_m$, is its minimal saturated neighborhood in $A_\a.$
\end{pf}

It is worth mentioning that due to Proposition \ref{trivialfacts} we
shall deal with entries $a_{ij}$ such that $0<a_{ij}<1.$

Let us present some criteria for the extremity of a matrix
$A_m\in{\mathbf{U}}_m.$

\begin{thm}\label{commoncriteria} Let $A_m$ be an element of ${\mathbf{U}}_m.$
Then the following assertions are equivalent:
\begin{itemize}
\item[(i)] The matrix $A_m$ is an extreme point of ${\mathbf{U}}_m;$
\item [(ii)] Every entry $ a_{ij}$ of $A_m$ with $0<a_{ij}<1$
has at least one saturated neighborhood, and minimal saturated neighborhoods of any two entries
$a_{ij},$ $a_{i'j'}$ with $0<a_{ij}, a_{i'j'}<1$ do not coincide.
\end{itemize}
\end{thm}

\begin{pf} (i)$\Rightarrow$ (ii). Let $A_m\in {\mathbf{Extr}\mathbf{U}}_m.$
Let us prove that every entry $a_{ij}$ of $A_m$ with $0<a_{ij}<1$
has at least one saturated neighborhood. We suppose the contrary,
i.e., there exist $0<a_{i_0j_0}<1$ having no saturated
neighborhoods, which means, for any $\a\in I_m$ with
$\a\supset\{i_0,j_0\}$ one has
\begin{eqnarray}\label{havenosaturated}
\sum\limits_{i,j\in \a}a_{ij}<|\a|.
\end{eqnarray}

We know that there are two cases either $i_0\neq j_0$ or $i_0=j_0.$
Let us consider the case $i_0\neq j_0,$ the second case can be
proceeded  by the same argument.

Since $I_m$ is a finite set, then a number of inequalities in
\eqref{havenosaturated} is finite. That is why, there exists
$0<\e_0<1$ such that
\begin{eqnarray*}
0<a_{i_0j_0}+\e_0<1, \, && \, 0<a_{i_0j_0}-\e_0<1,\\
\sum\limits_{i,j\in\a}a_{ij}+ 2\e_0<|\a|, \, && \, \sum\limits_{i,j\in\a}a_{ij}- 2\e_0<|\a|,
\end{eqnarray*}
for all $\a\supset\{i_0,j_0\}$.

Define two matrices $A'_m=(a'_{ij})_{i,j=1}^m$ and
$A''_m=(a''_{ij})_{i,j=1}^m,$ as follows
\begin{eqnarray*}
a_{ij}'&:=& \left\{\begin{array}{l}
                   a_{i_0j_0}+\e_0,  \quad(i,j)=(i_0,j_0)\\
                   a_{j_0i_0}+\e_0,  \quad(i,j)=(j_0,i_0)\\
                   a_{ij},  \quad\quad\quad (i,j)\neq(i_0,j_0)\vee (j_0,i_0)
                 \end{array}
\right.,\\[2mm]
a_{ij}''&:=& \left\{\begin{array}{l}
                   a_{i_0j_0}-\e_0, \quad (i,j)=(i_0,j_0)\\
                   a_{j_0i_0}-\e_0,  \quad(i,j)=(j_0,i_0)\\
                   a_{ij},  \quad\quad\quad (i,j)\neq(i_0,j_0)\vee (j_0,i_0)
                 \end{array}
\right..
\end{eqnarray*}
Due to the choice of $\e_0,$ we get $A'_m,A''_m\in{\mathbf{U}}_m$
and $2A_m=A'_m+A''_m$ which refutes the extremity of $A_m.$

Let us show that, for any two entries $a_{ij},$ $a_{i'j'}$ with
$0<a_{ij}, a_{i'j'}<1$, their corresponding minimal saturated
neighborhoods do not coincide. We suppose the contrary, i.e.,  for
two entries $0<a_{i_0j_0}<1$ and $0<a_{i'_0j'_0}<1$ their
corresponding minimal saturated neighborhoods coincide. This means
that there exists $\a_0\subset I_m$ with
$\{i_0,j_0,i'_0,j'_0\}\subset\a_0$ such that
\begin{eqnarray*}
\sum\limits_{i,j\in\a_0}a_{ij}=|\alpha_0|.
\end{eqnarray*}

We suppose that $i_0\neq j_0,$ $i'_0\neq j'_0.$ For the other
possible cases one can use the similar argument.

Since $A_{\a_0}$ is a common minimal saturated neighborhood of
entries $a_{i_0j_0}$ and $a_{i'_0j'_0}$, then due to the minimality of $\a_0$ it is clear that, for any $\a$ with
either $\a\supset\{i_0,j_0\}$ or $\a\supset\{i'_0,j'_0\}$ and
\begin{eqnarray*}
\sum\limits_{i,j\in\a}a_{ij}=|\a|,
\end{eqnarray*}
we have $\a_0\subset\a.$ In other words,
any saturated neighborhoods of $a_{i_0j_0}$ or $a_{i'_0j'_0}$ contain both of them.

Let us consider all index sets
$\beta\supset\{i_0,j_0\}$ and $\beta'\supset\{i'_0,j'_0\}$ such that
\begin{eqnarray}
\label{alphaiojo}\sum\limits_{i,j\in\beta}a_{ij}<|\beta|,\ \ \ \
\sum\limits_{i,j\in\beta'}a_{ij}<|\beta'|.
\end{eqnarray}
Since $I_m$ is a finite set then a number of inequalities in
\eqref{alphaiojo} is finite. Therefore, one can find $0<\e_0<1$ such
that
\begin{eqnarray*}
0<a_{i_0j_0}\pm\e_0<1, \, && \, 0<a_{i'_0j'_0}\pm\e_0<1,\\
\sum\limits_{i,j\in\beta}a_{ij}\pm 2\e_0<|\beta|, \, && \, \sum\limits_{i,j\in\beta'}a_{ij}\pm2\e_0<|\beta'|,
\end{eqnarray*}
for all $\beta\supset\{i_0,j_0\}$ and $\beta'\supset\{i'_0,j'_0\}$
satisfying inequalities \eqref{alphaiojo}.

Define two matrices $A'_m=(a'_{ij})_{i,j=1}^m$ and
$A''_m=(a''_{ij})_{i,j=1}^m,$ as follows
\begin{eqnarray*}
a_{ij}'&:=& \left\{\begin{array}{l}
                   a_{i_0j_0}+\e_0,  \quad(i,j)=(i_0,j_0)\\
                   a_{j_0i_0}+\e_0,  \quad(i,j)=(j_0,i_0)\\
                   a_{i'_0j'_0}-\e_0,  \quad(i,j)=(i'_0,j'_0)\\
                   a_{j'_0i'_0}-\e_0,  \quad(i,j)=(j'_0,i'_0)\\
                   a_{ij},  \ \ \ \quad\quad\quad (i,j)\neq(i_0,j_0)\vee (j_0,i_0)\vee (i'_0,j'_0)\vee(j'_0,i'_0)
                 \end{array}
\right.\\[2mm]
a_{ij}''&:=& \left\{\begin{array}{l}
                   a_{i_0j_0}-\e_0,  \quad (i,j)=(i_0,j_0)\\
                   a_{j_0i_0}-\e_0,  \quad(i,j)=(j_0,i_0)\\
                   a_{i'_0j'_0}+\e_0,  \quad (i,j)=(i'_0,j'_0)\\
                   a_{j'_0i'_0}+\e_0,  \quad(i,j)=(j'_0,i'_0)\\
                   a_{ij},  \ \ \ \quad\quad\quad (i,j)\neq(i_0,j_0)\vee (j_0,i_0)\vee (i'_0,j'_0)\vee(j'_0,i'_0)
                 \end{array}
\right..
\end{eqnarray*}
Let us show that $A'_m,A''_m\in{\mathbf{U}}_m.$ In fact, the
matrices $A'_m,A''_m$ are symmetric. We shall show that
$A'_m\in{\mathbf{U}}_m.$ By the same argument, one can show that
$A''_m\in{\mathbf{U}}_m.$

Let $\beta$ be any subset of $I_m.$ We want to estimate the sum
$\sum\limits_{i,j\in\beta}a'_{ij}.$ We are going to consider the
following cases.

\textsc{Case I.} Let $\beta\supseteq\{i_0,j_0\}$ and
$\beta\supseteq\{i'_0,j'_0\}.$ Then one gets
$\sum\limits_{i,j\in\beta}a'_{ij}=\sum\limits_{i,j\in\beta}a_{ij}\leq|\beta|.$

\textsc{Case II.} Let $\beta\nsupseteq\{i_0,j_0\}$ and
$\beta\nsupseteq\{i'_0,j'_0\}.$ Then
$\sum\limits_{i,j\in\beta}a'_{ij}=\sum\limits_{i,j\in\beta}a_{ij}\leq|\beta|.$

\textsc{Case III.} Let $\beta\nsupseteq\{i_0,j_0\}$ and
$\beta\supseteq\{i'_0,j'_0\}.$ Then
$\sum\limits_{i,j\in\beta}a'_{ij}=\sum\limits_{i,j\in\beta}a_{ij}-2\e_0<|\beta|.$

\textsc{Case IV.} Let $\beta\supseteq\{i_0,j_0\}$ and
$\beta\nsupseteq\{i'_0,j'_0\}.$ Then t
$\sum\limits_{i,j\in\beta}a'_{ij}=\sum\limits_{i,j\in\beta}a_{ij}+2\e_0.$
It is clear that the set $\beta$ could not be a saturated index in
the matrix $A_m$. In fact, every saturated index in the matrix $A_m$
containing the set $\{i_0,j_0\}$ should contain the set
$\{i'_0,j'_0\}.$ However, $\beta\nsupseteq\{i'_0,j'_0\}.$ Therefore,
due to choice of $\e_0$ we have
$\sum\limits_{i,j\in\beta}a'_{ij}=\sum\limits_{i,j\in\beta}a_{ij}+2\e_0<|\beta|.$

Consequently, the constructed matrices $A'_m,A''_m$ belong to
${\mathbf{U}}_m,$ and $2A_m=A'_m+A''_m$ which refutes the extremity
of $A_m.$

(ii)$ \Rightarrow $ (i). Suppose that the assertions (ii) are
satisfied. Let us prove $A_m\in{\mathbf{Extr}\mathbf{U}}_m$. We are
going to show it by using induction with respect to the order of the
matrix $A_m.$

Let $m=2.$ If the entries of the matrix $A_2\in{\mathbf{U}}_2$ are
either 1 or 0,
 then according to Proposition \ref{trivialfacts}
(iii) $A_2$ is extreme. Suppose that there exist at least one entry
$a_{ij}$ of $A_2$ with $0<a_{ij}<1.$  Then it is obvious that such
kind of entries' saturated neighborhood's radius is greater or equal
to 2. Consequently, $A_2$ is a minimal saturated neighborhood for
all $0<a_{ij}<1.$ It follows from (ii) that there is only one entry
$0<a_{ij}<1$ (of course, we only consider entries $a_{ij}$ with
$i\le j$) and the rest entries are either 1 or 0. After small
algebraic manipulations, we make sure that $A_2$  is an extreme
matrix in ${\mathbf{U}}_2.$

We suppose that the assumption of the theorem is true for all $m\le
k-1$, and we prove it for $m=k.$

If the entries of the matrix $A_k\in{\mathbf{U}}_k$ are either 1 or
0 then according to Proposition \ref{trivialfacts} (iii) $A_k$ is
extreme. Suppose that there exist some entries $0<a_{ij}<1.$  It
follows from (i) that every entry $0<a_{ij}<1$ has a minimal
saturated neighborhood, and we denote it by $A_{\a(a_{ij})}$, its
radius by $r(a_{ij}).$

Let us consider entries $a_{ij}$ with $0<a_{ij}<1,$ $r(a_{ij})\le k-1.$ Since
$A_{\a(a_{ij})}$ is a saturated principal submatrix of $A_m,$ then
according to Lemma \ref{indepenofminimsaturneigh}, the minimal
saturated neighborhoods $A_{\a(a_{i'j'})}$ of entries $0<a_{i'j'}<1$
in $A_{\a(a_{ij})}$  are  principal submatrices of $A_{\a(a_{ij})}.$
So, theorem assertions (i),(ii) are satisfied, for $A_{\a(a_{ij})}$,
and order of such a  matrix $A_{\a(a_{ij})}$ is less or equal to
$k-1$. Hence, according to the assumption of the induction we obtain
that $A_{\a(a_{ij})}\in {\mathbf{Extr}\mathbf{U}}_{|\a(a_{ij})|}.$
So, minimal saturated neighborhoods $A_{\a(a_{ij})}$ of entries
$0<a_{ij}<1$ of $A_k$ with $r(a_{ij})\le k-1$  are extreme in
${\mathbf{U}}_{|\a(a_{ij})|}.$

Let us consider  entries $a_{ij}$ with $0<a_{ij}<1,$ $r(a_{ij})= k.$ It
follows from (ii) that there is only one such kind of entries.

Suppose that the matrix can be decomposed as $2A_k=A_k'+A_k'',$
where $A'_k,A''_k\in{\mathbf{U}}_k.$ If the entries $a_{ij}$ of
$A_k$ are either 1 or 0 then according to Proposition
\ref{trivialfacts} $\rm(ii)$ we get $a_{ij}=a_{ij}'=a_{ij}''.$
Further, for the entries $0<a_{ij}<1$ with $r(a_{ij})\le k-1,$ since
their minimal saturated neighborhoods $A_{\a(a_{ij})}$ are extreme
in ${\mathbf{U}}_{|\a(a_{ij})|},$ according to Lemma
\ref{extremeprinsubmatix} we have
$A_{\a(a_{ij})}=A'_{\a(a_{ij})}=A''_{\a(a_{ij})},$ particularly,
$a_{ij}=a_{ij}'=a_{ij}''.$ Now, we must to show
$a_{ij}=a_{ij}'=a_{ij}''$ for an entry $0<a_{ij}<1$ with $r(a_{ij})=
k.$ We already mentioned that there is only one such kind of
entries, we denote it by $a_{i_0j_0}$ and its minimal saturated
neighborhood is $A_k.$ Since $2A_k=A_k'+A_k''$ then $A'_k$ and
$A''_k$ are also saturated matrices. We already know that for
entries $a_{ij}$ of $A_k$ with $(i,j)\neq(i_0,j_0)$ and
$(i,j)\neq(j_0,i_0),$  one has $a_{ij}=a_{ij}'=a_{ij}''.$ Therefore,
we get
\begin{eqnarray*}
&&0=k-k=\sum\limits_{i,j=1}^ka_{ij}-\sum\limits_{i,j=1}^ka'_{ij}=2(a_{i_0j_0}-a'_{i_0j_0}),\\
&&0=k-k=\sum\limits_{i,j=1}^ka_{ij}-\sum\limits_{i,j=1}^ka''_{ij}=2(a_{i_0j_0}-a''_{i_0j_0}).
\end{eqnarray*}
Consequently,  $a_{i_0j_0}=a'_{i_0j_0}=a''_{i_0j_0}.$ All these
facts bring to a conclusion that $A_k$ is a extreme matrix in
${\mathbf{U}}_k.$
\end{pf}

\begin{thm}\label{localcriterion} Let $A_m$ be an element of ${\mathbf{U}}_m.$
Then the following conditions are equivalent
\begin{itemize}
\item[(i)] The matrix $A_m$ is an extreme point of ${\mathbf{U}}_m;$
\item [(ii)] Every entry $ a_{ij}$ of $A_m$ with $0<a_{ij}<1$  has at least
one saturated neighborhood and its minimal saturated
neighborhood $A_{\a}$ is extreme in ${\mathbf{U}}_{|\a|};$
\item [(iii)] Every entry $ a_{ij}$ of $A_m$ with $0<a_{ij}<1$
has at least one saturated neighborhood and any its saturated
principal submatrix $A_{\a}$ is extreme in
$\mathbf{U}_{|\a|}.$
\end{itemize}
\end{thm}

\begin{pf} The implication (iii)$\Rightarrow$ (ii) is obvious.
Consider the implication (ii)$\Rightarrow$ (i). Assume that $A_m$
has the following decomposition $2A_m=A_m'+A_m''$ with
$A_m',A_m''\in{\mathbf{U}}_m.$ We know that if $a_{ij}=1\vee 0$ then
Proposition \ref{trivialfacts} (ii) yields
$a_{ij}=a_{ij}'=a_{ij}''.$ If $0<a_{ij}<1,$ since its minimal
saturated neighborhood $A_{\a}$ is extreme in $\mathbf{U}_{|\a|}$,
then due to  Proposition \ref{extremeprinsubmatix} one gets
$A_{\a}=A_{\a}'=A_{\a}'',$ particularly, $a_{ij}=a_{ij}'=a_{ij}''.$
These mean $A_m=A'_m=A''_m,$ i.e. $A_m$ is a extreme matrix in
${\mathbf{U}}_m.$

(i)$\Rightarrow$ (ii).
Let $A_m\in {\mathbf{Extr}\mathbf{U}}_m.$ Then according to Theorem
\ref{commoncriteria} every entry $ a_{ij}$ with $0<a_{ij}<1$ of
$A_m$ has at least one saturated neighborhood. We want to show that
a minimal saturated neighborhood $A_{\a}$ of such an entry is
extreme in $\mathbf{U}_{|\a|}$. Since $A_{\a}$
is a saturated principal submatrix of $A_m$ then using Lemma
\ref{indepenofminimsaturneigh} we deduce that, the minimal saturated neighborhood in $A_m$ of an entry $a_{i'j'}$ of $A_{\a}$ with $0<a_{i'j'}<1$ coincides with its minimal saturated neighborhood in $A_{\a}.$
Therefore, by applying Theorem \ref{commoncriteria} for the matrix
$A_{\a}$ we conclude that $A_{\a}\in \mathbf{Extr}\mathbf{U}_{\a}.$

(i)$\Rightarrow$ (iii). Let $A_m\in {\mathbf{Extr}\mathbf{U}}_m.$
Then according to Theorem \ref{commoncriteria} every entry $ a_{ij}$
of $A_m$ with $0<a_{ij}<1$  has at least one saturated neighborhood.
We want to show that any its saturated principal submatrix $A_{\a}$
is extreme in $\mathbf{U}_{|\a|}.$ Let $A_{\beta}$ be a saturated
principal submatrix  of $A_m\in{\mathbf{U}}_m.$ Since $A_{\beta}$ is
a saturated principal submatrix then using Lemma
\ref{indepenofminimsaturneigh} we can conclude that, for every entry
$0<a_{i'j'}<1$ of $A_{\beta},$ its minimal saturated neighborhood in
$A_m$ coincides with its minimal saturated neighborhood in
$A_{\beta}.$ Therefore, if we apply (ii)$\Rightarrow $(i) to the
submatrix $A_\beta$ we get
$A_{\beta}\in{\mathbf{Extr}\mathbf{U}}_{\beta}.$
\end{pf}

We are going to describe all extreme points of ${\mathbf{U}}^m$. It
is clear that the set ${\mathbf{U}}^m$ is a set of all the saturated
matrices of the set ${\mathbf{U}}_m.$

\begin{prop}\label{extrUmandextrUm} A matrix $A_m\in{\mathbf{U}}^m$
is extreme in ${\mathbf{U}}^m$ if and only if  $A_m$ is  extreme in
${\mathbf{U}}_m.$ Namely, one has
\begin{eqnarray*}
{\mathbf{Extr}}{\mathbf{U}}^m={\mathbf{U}}^m\cap
{\mathbf{Extr}}{\mathbf{U}}_m.
\end{eqnarray*}
\end{prop}

The proof of this proposition is straightforward.

From Theorem \ref{localcriterion} and Proposition
\ref{extrUmandextrUm} we immediately get

\begin{thm}\label{localcriterionforU^m} Let $A_m$ be an element of ${\mathbf{U}}^m.$
Then the following conditions are equivalent
\begin{itemize}
\item[(i)] The matrix $A_m$ is an extreme point of ${\mathbf{U}}^m;$
\item [(ii)] Every minimal saturated neighborhood $A_{\a}$ of $A_m$ is extreme in ${\mathbf{U}}^{|\a|};$
\item [(iii)] Every saturated principal submatrix $A_{\a}$ of $A_m$ is extreme in
$\mathbf{U}^{|\a|}.$
\end{itemize}
\end{thm}

It is worth mentioning that in the paper \cite{GS} the part $(i)\Rightarrow (iii)$ of Theorem \ref{localcriterionforU^m} was stated without any justification.

\section{Some properties of extreme points of the sets ${\mathbf{U}}_m$ and ${\mathbf{U}}^m$}\label{section3}

In this section, we are going to study some properties of extreme
matrices of the sets ${\mathbf{U}}_m$ and ${\mathbf{U}}^m$. The
results of this section will be used to solve Birkhoff's problem for
quadratic $G-$doubly stochastic operators.

Let us introduce the following sets
\begin{eqnarray*}
{\mathbf{U}}_m^{(0,1)}&:=&\Bigl\{A_m\in{\mathbf{U}}_m |
\quad a_{ij}=0\vee 1 \quad \forall \, i,j\in I_m\Bigr\},\\
{\mathbf{U}}_m^{(0,\frac12)}&:=&\left\{A_m\in{\mathbf{U}}_m | \quad a_{ii}=0, \quad a_{ij}=0\vee \frac12 \quad \forall \, i,j\in I_m\right\},\\
{\mathbf{U}}_m^{(0, \frac12, 1)}&:=&\left\{A_m\in{\mathbf{U}}_m |
\quad a_{ii}=0\vee 1, \quad a_{ij}=0\vee \frac12\vee 1 \quad \forall
\, i,j\in I_m\right\}.
\end{eqnarray*}

\begin{rem}\label{suminteger} The following assertions are evident:
\begin{itemize}
  \item [(i)] ${\mathbf{U}}_m^{(0,1)}\subset {\mathbf{U}}_m^{(0, \frac12, 1)} $ and
  ${\mathbf{U}}_m^{(0,\frac12)} \subset{\mathbf{U}}_m^{(0, \frac12, 1)}$;
  \item [(ii)] If $A\in {\mathbf{U}}_m^{(0, \frac12, 1)}$
  then $\sum\limits_{i,j\in\a}a_{ij}$ is an integer, for any $\a\subset I_m.$
\end{itemize}
\end{rem}

\begin{thm}\label{extrsubU0121} If $A_m\in{\mathbf{Extr}\mathbf{U}}_m$
then $a_{ii}=0\vee 1$ and $a_{ij}=0\vee\frac12\vee 1$ for any $i\neq
j$ and $i,j\in I_m.$ In other words, we have $
{\mathbf{Extr}}{\mathbf{U}}_m\subset{\mathbf{U}}_m^{(0, \frac12,
1)}. $
\end{thm}

\begin{pf} It is enough to show that $a_{ii}$ and $2a_{ij}$ are integers,
for any $i\neq j$ and $i,j\in I_m.$ We prove it by using induction
with respect to the order of $A_m$.

Let $m=2.$ Since $A_2\in {\mathbf{Extr}\mathbf{U}}_2$, then
according to Theorem \ref{commoncriteria}, there is at most one
entry $a_{ij}$ with $0<a_{ij}<1.$ If there is no such an entry then
the claim is obvious. Assume that there exists an entry
$a_{ij}$ with $0<a_{ij}<1.$ Then the matrix $A_2$ should be saturated and from
$a_{11}+2a_{12}+a_{22}=2$ we deduce that $a_{11},2a_{12},a_{22}$ are
integers.

Now suppose that the assertion of the theorem is true, for all
matrices $A_m\in{\mathbf{Extr}\mathbf{U}}_m$ of order $m\le k-1.$ We
prove it for matrices $A_m\in{\mathbf{Extr}\mathbf{U}}_m$ of order
$m=k.$

Since $A_k\in {\mathbf{Extr}\mathbf{U}}_m$ then, due to Theorem
\ref{localcriterion}, every entry $0<a_{ij}<1$ of $A_k$ has a
minimal saturated neighborhood $A_\a$ which is extreme in
$\mathbf{U}_{|\a|}.$ So, according to the assumption of induction,
for those minimal saturated neighborhoods with radius less or equal
to $k-1$ their entries $a_{ii}$ and $2a_{ij}$ are integers, for any
$i\neq j$ and $i,j\in I_k.$ In other words, the assumption of
induction allow us to say that  all $a_{ii}$ and $2a_{ij}$ ($i<j$)
of entries of $A_k$ are integers except which has a minimal
saturated neighborhood with radius $k.$

Now, assume that there is an entry $a_{ij}$ with $0<a_{ij}<1$ which
has saturated neighborhoods with radius equal to $k.$ Then according
to Theorem \ref{commoncriteria}, there is only one such an entry, we
denote it by $a_{i_0j_0}$, and the matrix $A_k$ should be saturated.
In this case, we already know that for entries $a_{ij}$ of $A_k$
with $(i,j)\neq(i_0,j_0),$ if $i=j$ then $a_{ii}$  is an integer,
and if $i\neq j$ then $2a_{ij}$ is an integer.  Then from
\begin{eqnarray*}
\sum\limits_{i,j=1}^ka_{ij}=\sum\limits_{i=1}^ka_{ii}+2\sum\limits_{i<j}a_{ij}=k
\end{eqnarray*}
we conclude that if $i_0=j_0$ then $a_{i_0i_0}$ is an integer, and
if $i_0\neq j_0$ then $2a_{i_0j_0}$ is an integer as well. All of
these mean that $a_{ii}$ and $2a_{ij}$ are integers, for any $i\neq
j$ and $i,j\in I_k.$
\end{pf}

Now, combining Proposition \ref{trivialfacts} (iii) with Theorem
\ref{extrsubU0121} we have
\begin{eqnarray}
{\mathbf{U}}_m^{(0,1)}\subset {\mathbf{Extr}}{\mathbf{U}}_m\subset
{\mathbf{U}}_m^{(0, \frac12, 1)}.
\end{eqnarray}

\begin{rem}\label{Extrm=2} If $m=2$ then it is easy to show the following equalities:
\begin{itemize}
  \item [(i)] ${\mathbf{U}}_2^{(0,\frac12)}\cap {\mathbf{Extr}}{\mathbf{U}}_2=\circleddash_{2\times 2};$
  \item [(ii)] ${\mathbf{U}}_2^{(0,\frac12)}\cup {\mathbf{Extr}}{\mathbf{U}}_2
  ={\mathbf{U}}_2^{(0, \frac12, 1)};$
\end{itemize}
here, as before, $\circleddash_{2\times 2}$ is $2\times 2$ zero
matrix.
\end{rem}

Our next aim is that Remark \ref{Extrm=2} (i) holds true for any $m.$ To this end, we introduce the following
useful conception.

\begin{defin}
A matrix $A_m\in {\mathbf{U}}_m^{(0, \frac12, 1)}$ is said to be
$F_m-$matrix if $A_m$ is saturated, and having no other saturated
principal submatrices.
\end{defin}

\begin{prop}\label{FmmatrixandU012}
Let $A_m$ be an $F_m-$matrix. If $m\ge 3$ then $A_m\in
{\mathbf{U}}_m^{(0,\frac12)}.$
\end{prop}

\begin{pf} Suppose that $A_m$ is an $F_m-$matrix
and $m\ge 3.$ If $\a=\{i\}$ then it follows from the definition of
the $F_m-$matrix and Remark \ref{suminteger} (ii) that $a_{ii}\le
0,$ which means $a_{ii}=0$ for all $i\in I_m.$ If $\a=\{i,j\}$ with
$i\neq j$, then with the same reason as the previous case, from the
inequality
$$\sum\limits_{i,j\in\alpha}a_{ij}=a_{ii}+a_{jj}+2a_{ij}\le 1,$$
we find out that $a_{ij}=0\vee \frac12$ for all $i\neq j$ and
$i,j\in I_m.$ This means that $A_m\in {\mathbf{U}}_m^{(0,\frac12)}.$
\end{pf}

\begin{rem}
For $m=2,$ there is only one $F_2-$matrix which is
$\left(
  \begin{array}{cc}
    0 & 1 \\
    1 & 0 \\
  \end{array}
\right).$
\end{rem}

From now on, we shall only consider the case $m\ge 3.$

\begin{exam} Let us consider the following $m\times m$ matrix
\begin{equation*}
N_m=\left(
      \begin{array}{cccccccc}
        0 & \frac12 & 0 & 0 & \cdots & 0 & 0 & \frac12 \\
        \frac12 & 0 & \frac12 & 0 & \cdots & 0  & 0 & 0 \\
        0 & \frac12 & 0 & \frac12 & \cdots & 0 & 0 & 0 \\
        \cdots & \cdots & \cdots & \cdots & \cdots & \cdots & \cdots & \cdots \\
        0 & 0 & 0 & 0 & \cdots & \frac12 & 0 & \frac12 \\
        \frac12 & 0 & 0 & 0 & \cdots & 0 & \frac12 & 0 \\
      \end{array}
    \right).
\end{equation*}
It is easy to see that if $m\ge 3$ then  $N_m$ is an $F_m-$matrix.
\end{exam}

\begin{prop}\label{Fmmatrixnotextreme}
Let $A_m$ be an $F_m-$matrix. If $m\ge 3$ then $A_m$ is not an
extreme element of ${\mathbf{U}}_m.$
\end{prop}

\begin{pf} Suppose that $A_m$ is an $F_m-$matrix and $m\ge 3$. Then, due to
Proposition \ref{FmmatrixandU012}, we have $A_m\in
{\mathbf{U}}_m^{(0,\frac12)}.$ By the definition of the
$F_m-$matrix, the number of entries of $A_m$, which are $1/2$, is
equal to $m$ (of course, we are only speaking about such entries
$a_{ij}$ in which $i\le j$). On the other hand, a minimal saturated
neighborhood of any nonzero entry is $A_m$ itself. So, there are at
least two nonzero entries whose minimal saturated neighborhoods
coincide. Then Theorem \ref{commoncriteria} brings to a conclusion
that the matrix $A_m$ is not extreme in ${\mathbf{U}}_m.$
\end{pf}

\begin{prop}
Let $A_m\in {\mathbf{U}}_m^{(0, \frac12, 1)}.$ If $A_m$ has a
principal $F_k-$submatrix, where $k\ge 3,$ then $A_m$ is not extreme
in ${\mathbf{U}}_m$.
\end{prop}

\begin{pf} Let $A_m\in{\mathbf{U}}_m^{(0, \frac12, 1)}$ and $A_k$ be a principal
$F_k-$submatrix of $A_m.$ Since $k\ge 3$ then, due to Propositions
\ref{FmmatrixandU012} and \ref{Fmmatrixnotextreme}, we get
$A_k\in{\mathbf{U}}_k^{(0, \frac12)}$ and
$A_k\notin{\mathbf{Extr}}{\mathbf{U}}_k.$ Then a number of entries
of $A_k$, which are equal to $\dfrac12$, is equal to $k$ and a
minimal saturated neighborhood of such an entry is $A_k$, which is
not extreme in ${\mathbf{U}}_k.$ Then, according to Theorem
\ref{localcriterion}, we conclude that the matrix $A_m$ is not
extreme in ${\mathbf{U}}_m.$
\end{pf}

Now, we are ready to formulate one of important properties of the
extremal matrices in ${\mathbf{U}}_m.$

\begin{thm}\label{Extrforanym}  The following equality is satisfied for any $m$
\begin{eqnarray}\label{interExtrUmU012}
{\mathbf{Extr}}{\mathbf{U}}_m\cap
{\mathbf{U}}_m^{(0,\frac12)}=\circleddash_{m\times m},
\end{eqnarray}
here, as before, $\circleddash_{m\times m}$ is zero matrix.
\end{thm}
\begin{pf}
Let $m\ge 3,$ otherwise Remark \ref{Extrm=2} $\rm(i)$ yields the
assertion.  Suppose that $A_m\in {\mathbf{Extr}}{\mathbf{U}}_m$ and
$A_m\neq \circleddash_m.$  We want to show that $A_m\notin
{\mathbf{U}}_m^{(0,\frac12)}.$ We assume that
$A_m\in{\mathbf{U}}_m^{(0,\frac12)}.$ Since $A_m$ is a nonzero
extreme matrix in ${\mathbf{U}}_m^{(0,\frac12)}$, then due to
Theorem \ref{localcriterion}, there exists at least one entry
$a_{ij}=\frac12$ and its minimal saturated neighborhood $A_{\a}$
should be extreme in $\mathbf{U}_{\a}.$ On the other hand, since
$A_m\in{\mathbf{U}}_m^{(0,\frac12)}$ then the minimal saturated
neighborhood $A_{\a}$ of  that entry (i.e. $a_{ij}=\frac12$) is
$F_{|\a|}-$matrix which is not extreme $\mathbf{U}_{\a}.$ This
contradiction proves the desired assertion.
\end{pf}

\begin{cor}\label{bestcriterion} Let $A_m\in{\mathbf{Extr}}{\mathbf{U}}_m$,
then the following assertions hold true:
\begin{enumerate}
  \item [(i)] One has $A_m\in {\mathbf{U}}_m^{(0, \frac12, 1)};$
  \item [(ii)] Every entry $a_{ij}=\frac12$ has at least one saturated neighborhood;
  \item [(iii)] For any saturated principal submatrix
$A_{\a}$, one has $A_{\a}\notin{\mathbf{U}}_{|\a|}^{(0,\frac12)}$.
\end{enumerate}
\end{cor}

By using Corollary \ref{bestcriterion} and Proposition
\ref{extrUmandextrUm}, we get the following description of the
extreme points of ${\mathbf{U}}^m.$

\begin{cor}\label{bestcriterionforU^m} Let $A_m\in{\mathbf{Extr}}{\mathbf{U}}^m$.
Then the following assertions hold true:
\begin{enumerate}
  \item [(i)]  One has $A_m\in {\mathbf{U}}_m^{(0, \frac12, 1)};$
  \item [(ii)] For any saturated principal submatrix $A_{\a}$, one has
$A_{\a}\notin{\mathbf{U}}_{|\a|}^{(0,\frac12)}$.
\end{enumerate}
\end{cor}

It seems the following conjectures hold true.

\begin{conjec}
Let $A_m\in{\mathbf{Extr}}{\mathbf{U}}_m$ if and only if the
following assertions hold true:
\begin{enumerate}
  \item [(i)] One has $A_m\in {\mathbf{U}}_m^{(0, \frac12, 1)};$
  \item [(ii)] Every entry $a_{ij}=\frac12$ has at least one saturated neighborhood;
  \item [(iii)] For any saturated principal submatrix $A_{\a}$ one has
  $A_{\a}\notin{\mathbf{U}}_{|\a|}^{(0,\frac12)}$.
\end{enumerate}
\end{conjec}

\begin{conjec}
Let $A_m\in{\mathbf{U}}^m$. Then $A_m\in{\mathbf{Extr}}{\mathbf{U}}^m$
if and only if the following assertions hold true:
\begin{enumerate}
  \item [(i)]  One has $A_m\in {\mathbf{U}}_m^{(0, \frac12, 1)};$
  \item [(ii)] For any saturated principal submatrix $A_{\a}$, one has
$A_{\a}\notin{\mathbf{U}}_{|\a|}^{(0,\frac12)}$.
\end{enumerate}
\end{conjec}

\section{Canonical forms of extreme points of the sets ${\mathbf{U}}_m$ and ${\mathbf{U}}^m$ }\label{section4}

In this section, we are going to describe location of nonzero
entries of extreme matrices and canonical forms of extreme matrices
of the sets ${\mathbf{U}}_m$ and ${\mathbf{U}}^m$. Based on the
canonical forms of extreme points we are going to study an algebraic
structure of the sets ${\mathbf{U}}_m$ and ${\mathbf{U}}^m$ (see
sec. \ref{section5}).

\begin{prop} Let $A_m\in{\mathbf{Extr}}{\mathbf{U}}_m$. Then the following statements hold true:
\begin{itemize}
  \item [(i)]\label{minimalsatur} Every nonzero entry
of $A_m$ has a minimal saturated neighborhood;
  \item [(ii)]\label{maximalsatur} Every nonzero entry
of $A_m$ has a unique common maximal saturated neighborhood.
\end{itemize}
\end{prop}

\begin{pf} $\rm(i).$ Due to Corollary
\ref{bestcriterion} (i) and (ii), we have that $A_m\in
{\mathbf{U}}_m^{(0, \frac12, 1)}$ and every entry $a_{ij}=\dfrac12$
has a minimal saturated neighborhood. Now, let us consider such
entries with $a_{ij}=1.$ Then, there are two cases either $i=j$ or
$i\neq j.$ In both cases, a submatrix $A_\a$ with $\a=\{i,j\}$ of
$A_m$ is a minimal saturated neighborhood of $a_{ij}=1$.

$\rm(ii).$ Due to (i) every nonzero entry of $A_m$ has a minimal
saturated neighborhood. By $\a_1,\a_2,\cdots,\a_k$ we denote
saturated index sets, corresponding to these minimal saturated
neighborhoods of nonzero entries of $A_m$. According to Proposition
\ref{unionandinterofsatutset} (ii), an index set
$\a=\bigcup\limits_{i=1}^k\a_i$ is saturated. If we consider a
principal submatrix $A_\a$ of $A_m,$ corresponding to the saturated
index set $\a,$ then $A_\a$ is saturated. Since $A_\a$ contains all
nonzero entries of $A_m,$ therefore, $A_\a$ is a common maximal
saturated neighborhood of every nonzero entry of $A_m.$ The
uniqueness is trivial.
\end{pf}

\begin{cor}\label{alphaandalpha'}
 If $A_m\in{\mathbf{Extr}}{\mathbf{U}}_m$, then there exist two index sets
$\a$ and $\a'$ such that $\a\cup\a'=I_m,$ $\a\cap\a'=\emptyset,$ and
satisfying the following conditions:
\begin{itemize}
  \item [(i)] $A_\a$ is a saturated principal submatrix  of $A_m$;
  \item [(ii)] $A_\a$ contains all nonzero entries of $A_m;$
  \item [(iii)] $A_{\a'}=\circleddash_{|\a'|}$ and  $a_{ij'}=a_{i'j}=0$ for any $i,j\in\a,$ and $i',j'\in\a'.$
\end{itemize}
\end{cor}

\begin{lem}\label{criterionforsaturated}
Let $A_m\in{\mathbf{Extr}}{\mathbf{U}}_m$. A matrix $A_m$ is
saturated if and only if every row of $A_m$ has at least one nonzero
entry.
\end{lem}

\begin{pf} \emph{If part.} Let
$A_m\in{\mathbf{Extr}}{\mathbf{U}}_m$, and its every row has at
least one nonzero entry. Then, due to Corollary \ref{alphaandalpha'}
(i) and (ii), there exist two index sets $\a$ and $\a'$ with
$\a\cup\a'=I_m,$ $\a\cap\a'=\emptyset$ such that $A_\a$ is a
saturated principal submatrix  of $A_m$ containing all nonzero
entries of $A_m.$ Moreover, since every row of $A_m$ has at least
one nonzero entry, it follows form  Corollary \ref{alphaandalpha'}
(iii) that $\a'=\emptyset,$ which means $A_m=A_\a$ is a saturated
matrix.

\emph{Only if part.} Let $A_m$ be a saturated extreme matrix in
${\mathbf{U}}_m.$ We must to show that every row of $A_m$ has at
least one nonzero entry. Assume the contrary i.e. there is a row of
$A_m$ with zero entries. We denote it by $i_0.$ In this case,
symmetricity of $A_m$ implies that $i_0^{\rm {th}}$ column of $A_m$
has zero entries as well. Let $\a_0=I_m\setminus\{i_0\},$ then from
 $A_m\in\mathbf{U}_m,$ one gets
\begin{eqnarray*}
\sum\limits_{i,j\in\a_0}a_{ij}\le|\a_0|=m-1.
\end{eqnarray*}
On the other hand, since $A_m$ is saturated, it follows that
\begin{eqnarray*}
\sum\limits_{i,j=1}^ma_{ij}=\sum\limits_{i,j\in\a_0}a_{ij}+a_{i_0i_0}+2\sum\limits_{j\in\a_0}a_{i_0j}=\sum\limits_{i,j\in\a_0}a_{ij}=m.
\end{eqnarray*}
This contradiction shows that every row of $A_m$ has at least one nonzero entry.
\end{pf}

By using Proposition \ref{extrUmandextrUm} and Lemma
\ref{criterionforsaturated} one can get the following
\begin{cor}
Let $A_m\in{\mathbf{Extr}}{\mathbf{U}}_m$. Then the matrix $A_m$ is
extreme in ${\mathbf{U}}^m$ if and only if its every row has at
least one nonzero entry.
\end{cor}

\begin{cor}\label{nonsaturatedcase}
Let $A_m$ be a nonsaturated  extreme matrix in ${\mathbf{U}}_m$ and
$K:=m-\sum\limits_{i,j=1}^ma_{ij}.$ Then, there exist two index set
$\a$ and $\a'$ with $|\a'|=K,$ $\a=I_m\setminus\a',$ satisfying the
following conditions:
\begin{itemize}
  \item [(i)] $A_\a$ is a saturated principal submatrix of $A_m$;
  \item [(ii)] $A_\a$ contains all nonzero entries of $A_m;$
  \item [(iii)] Every $i$ row and $i$ column of $A_m$ consists zeroes, for any $i\in\a';$
  \end{itemize}
\end{cor}

We will introduce the following notation. Let $A_m$ be a matrix and
$\pi$ be a permutation of the set $I_m=\{1,2,\cdots,m\}.$ Define a
matrix as follows
$$A_{\pi(m)}=(a{'}_{ij})_{i,j=1}^m,\quad a{'}_{ij}=a_{\pi(i)\pi(j)} \ \forall \ i,j=\overline{1,m}.$$

\begin{prop}\label{permutasionofExtr}
Let $A_m$ be a matrix and $\pi$ be a permutation of the set $I_m=\{1,2,\cdots,m\}.$
Then the following assertions hold true:
\begin{itemize}
  \item [(i)] if $A_m\in {\mathbf{U}}_m$ then $A_{\pi(m)}\in {\mathbf{U}}_m;$
  \item [(ii)] if $A_m\in {\mathbf{Extr}}{\mathbf{U}}_m$ then $A_{\pi(m)}
\in {\mathbf{Extr}}{\mathbf{U}}_m;$
  \item [(iii)] if $A_m\in {\mathbf{U}}^m$ then $A_{\pi(m)}\in {\mathbf{U}}^m;$
  \item [$(iv)$] if $A_m\in {\mathbf{Extr}\mathbf{U}}^m$ then
$A_{\pi(m)}\in {\mathbf{Extr}\mathbf{U}}^m.$
\end{itemize}
\end{prop}
\begin{pf} (i). Let $A_m\in {\mathbf{U}}_m.$  Due to $a_{ij}=a_{ji}, \ \forall
\ i,j=\overline{1,m}$ we have
$$a{'}_{ij}=a_{\pi(i)\pi(j)}=a_{\pi(j)\pi(i)}=a{'}_{ji}, \ \forall \ i,j=\overline{1,m}.$$
This means that $(A_{\pi(m)})^{t}=A_{\pi(m)}.$ Let $\alpha\subset
I_m$, and put $\beta=\pi(\alpha)$. It is clear that
$|\beta|=|\alpha|.$ From $A_m\in{\mathbf{U}}_m$ it follows that
$$\sum\limits_{i,j\in\alpha}a{'}_{ij}=\sum\limits_{i,j\in\alpha}a_{\pi(i)\pi(j)}=\sum\limits_{i,j\in\beta}a_{ij}\le|\beta|=|\alpha|,$$
which implies  $A_{\pi(m)}\in {\mathbf{U}}_m.$

(ii). Let $A_m\in {\mathbf{Extr}\mathbf{U}}_m.$ We want to show that
$A_{\pi(m)}\in {\mathbf{Extr}\mathbf{U}}_m.$ We suppose the
contrary, i.e., there exist matrices $A{'}_m, A{''}_m\in
{\mathbf{U}}_m$ such that $2A_{\pi(m)}=A{'}_{m}+A{''}_{m}$ and
$A_{\pi(m)}\neq A{'}_{m},$ $A_{\pi(m)}\neq A{''}_{m}.$ Let us
consider the matrices $A{'}_{\pi^{-1}(m)}, A{''}_{\pi^{-1}(m)}$ then
one has
$$2A_m=A{'}_{\pi^{-1}(m)}+A{''}_{\pi^{-1}(m)},\quad A_{\pi^{-1}(m)}\neq A_{m},
\quad A_{\pi^{-1}(m)}\neq A_{m}.$$ This contradicts to $A_m\in
{\mathbf{Extr}\mathbf{U}}_m.$

(iii). Let $A_m\in {\mathbf{U}}^m$, then from
$\sum\limits_{i,j=1}^ma_{ij}=m$ one finds
$$\sum\limits_{i,j=1}^ma'_{ij}=\sum\limits_{i,j=1}^ma_{\pi(i)\pi(j)}
=\sum\limits_{i,j=1}^ma_{ij}=m,$$ which means  $A_{\pi(m)}\in
{\mathbf{U}}^m.$

Proposition \ref{extrUmandextrUm} with (ii),(iii) yields the
assertion (iv).
\end{pf}

By means of Proposition \ref{permutasionofExtr} and Corollary
\ref{nonsaturatedcase} we are going to provide a canonical form of
extreme points of ${\mathbf{U}}_m.$

Let $A_m\in {\mathbf{U}}_m$ be an extreme matrix in ${\mathbf{U}}_m$
and $k=\sum\limits_{i,j=1}^ma_{ij}.$ Then there exists a permutation
$\pi$ of the set $I_m$ such that the matrix $A_{\pi(m)}$ has the
following form
\begin{eqnarray}\label{canonicalformExtr}
A_{\pi(m)}=\left(
      \begin{array}{cc}
        A_k & \circleddash_{k\times m-k} \\
        \circleddash_{m-k\times k} & \circleddash_{m-k\times m-k} \\
      \end{array}
    \right),
\end{eqnarray}
here, as before, $\circleddash_{k\times m-k},$
$\circleddash_{m-k\times k},$ and $\circleddash_{m-k\times m-k}$ are
zero matrices and $A_k$ is an extreme saturated matrix in
${\mathbf{U}}^k$, i.e., $A_k\in{\mathbf{Extr}\mathbf{U}}^k.$ The
form \eqref{canonicalformExtr} is called \emph{a canonical form} of
the extreme matrix $A_m.$

In the sequel, without loss of generality, we will assume that an
extreme matrix $A_m$ has a canonical form \eqref{canonicalformExtr}.

Let $\alpha$ and $\beta$ be two nonempty disjoint partitions of the
set $I_m,$ i.e., $\alpha\cap\beta=\emptyset$ and
$\alpha\cup\beta=I_m.$ Let $A_\alpha=(a_{ij})_{i,j\in\alpha}$ and
$B_\beta=(b_{ij})_{i,j\in\beta}$ be two matrices. Define a matrix
$C_{\alpha\cup\beta}=(c_{ij})_{\alpha\cup\beta}$ as follows:
\begin{eqnarray}\label{Calphabeta}
c_{ij}=\left\{
\begin{array}{l}
  a_{ij} \quad i,j\in\alpha \\
  b_{ij} \quad i,j\in\beta \\
  0 \quad i\in\alpha, \ j\in\beta \\
  0 \quad i\in\beta, \ j\in\alpha.
\end{array}
\right.
\end{eqnarray}
\begin{prop}\label{AalphaBbetaCalphabeta}
Let $\alpha,$ $\beta$ be two nonempty disjoint partitions of $I_m$
and $A_\alpha=(a_{ij})_{i,j\in\alpha},$
$B_\beta=(b_{ij})_{i,j\in\beta}$ be two matrices. Let
$C_{\alpha\cup\beta}=(c_{ij})_{\alpha\cup\beta}$ be a matrix defined
by \eqref{Calphabeta}. Then the following assertions hold true:
\begin{itemize}
  \item[$(i)$] if $A_\alpha\in\mathbf{U}_{|\alpha|}$ and $B_\beta\in\mathbf{U}_{|\beta|}$ then $C_{\alpha\cup\beta}\in \mathbf{U}_{|\alpha|+|\beta|};$
  \item[$(ii)$] if $A_\alpha\in\mathbf{Extr}\mathbf{U}_{|\alpha|}$ and $B_\beta\in\mathbf{Extr}\mathbf{U}_{|\beta|}$ then $C_{\alpha\cup\beta}\in\mathbf{Extr}\mathbf{U}_{|\alpha|+|\beta|};$
  \item[$(iii)$] if $A_\alpha\in\mathbf{U}^{|\alpha|}$ and $B_\beta\in\mathbf{U}^{|\beta|}$ then $C_{\alpha\cup\beta}\in \mathbf{U}^{|\alpha|+|\beta|};$
  \item[$(iv)$] if $A_\alpha\in\mathbf{Extr}\mathbf{U}^{|\alpha|}$ and $B_\beta\in\mathbf{Extr}\mathbf{U}^{|\beta|}$ then $C_{\alpha\cup\beta}\in\mathbf{Extr}\mathbf{U}^{|\alpha|+|\beta|}.$
\end{itemize}
\end{prop}
\begin{pf} According to Proposition \ref{permutasionofExtr} we may
assume that $\alpha=\{1,2,\cdots,i\}$ and
$\beta=\{i+1,i+2,\cdots,m\}.$ Then the matrix
$C_{\alpha\cup\beta}=(c_{ij})_{\alpha\cup\beta}$ given by
\eqref{Calphabeta} has the following form
\begin{eqnarray}\label{canonicalCalphabeta}
C_{\alpha\cup\beta}=\left(
      \begin{array}{cc}
        A_\alpha & \circleddash_{|\alpha|\times |\beta|} \\
        \circleddash_{|\beta|\times |\alpha|} & B_\beta \\
      \end{array}
    \right).
\end{eqnarray}

(i). Let $A_\alpha\in\mathbf{U}_{|\alpha|}$ and
$B_\beta\in\mathbf{U}_{|\beta|}.$ We want to show that
$C_{\alpha\cup\beta}\in \mathbf{U}_{|\alpha|+|\beta|}.$ Indeed, it
follows from \eqref{canonicalCalphabeta} that
$(C_{\alpha\cup\beta})^{t}=C_{\alpha\cup\beta}.$ Let $\gamma\subset
I_m(=\alpha\cup\beta)$ be any subset of $I_m.$ Let
$\gamma_\alpha=\gamma\cap\alpha$, $\gamma_\beta=\gamma\cap\beta$
then $\gamma_\alpha\cap\gamma_\beta=\emptyset$,
$\gamma_\alpha\cup\gamma_\beta=\gamma.$ It is clear that
\begin{eqnarray*}
\sum\limits_{i,j\in\gamma}c_{ij}&=&\sum\limits_{i,j\in\gamma_\alpha}c_{ij}
+\sum\limits_{i,j\in\gamma_\beta}c_{ij}+\sum\limits_{i\in\gamma_\alpha, j\in\gamma_\beta}c_{ij}+
\sum\limits_{i\in\gamma_\beta,j\in\gamma_\alpha}c_{ij}\\
&=&\sum\limits_{i,j\in\gamma_\alpha}a_{ij}+\sum\limits_{i,j\in\gamma_\beta}b_{ij}\leq|\gamma_\alpha|+|\gamma_\beta|=|\gamma|
\end{eqnarray*}
This means that $C_{\alpha\cup\beta}\in \mathbf{U}_{|\alpha|+|\beta|}.$

$(ii).$ Let $A_\alpha\in\mathbf{Extr}\mathbf{U}_{|\alpha|}$ and
$B_\beta\in\mathbf{Extr}\mathbf{U}_{|\beta|}.$ We suppose that there
exist $C^{'}_{\alpha\cup\beta}$ and $C^{''}_{\alpha\cup\beta}$ such
that
$2C_{\alpha\cup\beta}=C^{'}_{\alpha\cup\beta}+C^{''}_{\alpha\cup\beta}.$
We may assume that the matrices $C^{'}_{\alpha\cup\beta}$ and
$C^{''}_{\alpha\cup\beta}$ have the following form
\begin{eqnarray*}
C^{'}_{\alpha\cup\beta}=\left(
      \begin{array}{cc}
        C^{'}_\alpha & C^{'}_{|\alpha|\times |\beta|} \\
        C^{'}_{|\beta|\times |\alpha|} & C^{'}_\beta \\
      \end{array}
    \right), \quad
C^{''}_{\alpha\cup\beta}=\left(
      \begin{array}{cc}
        C^{''}_\alpha & C^{''}_{|\alpha|\times |\beta|} \\
        C^{''}_{|\beta|\times |\alpha|} & C^{''}_\beta \\
      \end{array}
    \right),
\end{eqnarray*}
then it follows from \eqref{canonicalCalphabeta} that
\begin{eqnarray*}
2\circleddash_{|\alpha|\times |\beta|}&=&C^{'}_{|\alpha|\times |\beta|}+C^{''}_{|\alpha|\times |\beta|},\\
2\circleddash_{|\beta|\times |\alpha|}&=&C^{'}_{|\beta|\times |\alpha|}+C^{''}_{|\beta|\times |\alpha|},\\
2A_\alpha&=&C^{'}_\alpha+C^{''}_\alpha, \\
2B_\beta&=&C^{'}_\beta+C^{''}_\beta.
\end{eqnarray*}
Hence, we obtain $C^{'}_{|\alpha|\times
|\beta|}=C^{''}_{|\alpha|\times
|\beta|}=\circleddash_{|\alpha|\times |\beta|}$ and
$C^{'}_{|\beta|\times |\alpha|}=C^{''}_{|\beta|\times
|\alpha|}=\circleddash_{|\beta|\times |\alpha|}.$ From
$A_\alpha\in\mathbf{Extr}\mathbf{U}_{|\alpha|}$,
$B_\beta\in\mathbf{Extr}\mathbf{U}_{|\beta|}$ we find
$C^{'}_\alpha=C^{''}_\alpha=A_\alpha,$
$C^{'}_\beta=C^{''}_\beta=B_\beta.$ This means that
$C^{'}_{\alpha\cup\beta}=C^{''}_{\alpha\cup\beta}=C_{\alpha\cup\beta},$
i.e,
$C_{\alpha\cup\beta}\in\mathbf{Extr}\mathbf{U}_{|\alpha|+|\beta|}.$

(iii). Let $A_\alpha\in\mathbf{U}^{|\alpha|}$ and
$B_\beta\in\mathbf{U}^{|\beta|}.$ One can see that
$$\sum\limits_{i,j=1}^mc_{ij}=\sum\limits_{i,j\in\alpha}c_{ij}
+\sum\limits_{i,j\in\beta}c_{ij}=\sum\limits_{i,j\in\alpha}a_{ij}
+\sum\limits_{i,j\in\beta}b_{ij}=|\alpha|+|\beta|=m$$
this means that $C_{\alpha\cup\beta}\in \mathbf{U}^{|\alpha|+|\beta|}.$

The assertion (iv) immediately  follows from Proposition
\ref{extrUmandextrUm} and assertions (ii) and (iii).
\end{pf}

Now we are going to consider an extension problem: \emph{let
$A_m\in{\mathbf{U}}_m$ be a non-saturated matrix. Is there a
saturated matrix $A_{m+1}\in{\mathbf{U}}^{m+1}$ containing a matrix
$A_m$ as a principal sub-matrix? In other words, is it possible to
make a non-saturated matrix as a saturated matrix by increasing its
order?} If the extension problem has a positive answer, for the sets
${\mathbf{U}}_m$ and ${\mathbf{U}}^{m+1}$, then we use the following
natation ${\mathbf{U}}_m\hookrightarrow{\mathbf{U}}^{m+1}.$ We shall
solve this extension problem in a general setting.

\begin{prop}
Let $A_m\in{\mathbf{U}}_m$. Then there exists a saturated matrix
$A_{m+1}\in{\mathbf{U}}^{m+1}$ containing a matrix $A_m$ as a
principal sub-matrix, i.e.,
${\mathbf{U}}_m\hookrightarrow{\mathbf{U}}^{m+1}.$
\end{prop}

\begin{pf}
We shall prove the assertion in two steps.

\textsc{Step I.} Let us prove that if
$A_m\in{\mathbf{Extr}\mathbf{U}}_m$ then there exists a saturated
matrix $A_{m+1}\in{\mathbf{Extr}\mathbf{U}}^{m+1}$ containing a
matrix $A_m$ as a principal sub-matrix, i.e.,
${\mathbf{Extr}\mathbf{U}}_m\hookrightarrow{\mathbf{Extr}\mathbf{U}}^{m+1}.$
Indeed, we suppose that the matrix
$A_m\in{\mathbf{Extr}\mathbf{U}}_m$ has the following form
\begin{eqnarray*}
A_m=\left(
      \begin{array}{cc}
        A_k & \circleddash_{k\times m-k} \\
        \circleddash_{m-k\times k} & \circleddash_{m-k\times m-k} \\
      \end{array}
    \right),
\end{eqnarray*}
where $A_k\in{\mathbf{Extr}\mathbf{U}}^k$ and
$\sum\limits_{i,j=1}^ma_{ij}=k\leq m.$ Let us consider the following
matrix
\begin{eqnarray*}
A_{m+1-k}=\left(
      \begin{array}{cc}
        \circleddash_{m-k\times m-k} & \left(\frac12\right)_{m-k\times 1} \\
        \left(\frac12\right)_{1\times m-k} & 1 \\
      \end{array}
    \right),
\end{eqnarray*}
where $\left(\frac12\right)_{1\times
m-k}=\Bigl(\underbrace{\frac12,\cdots,\frac12}_{m-k}\Bigr)$ and
$\left(\frac12\right)_{m-k\times 1}=\left(\frac12\right)_{1\times
m-k}^t.$ It is clear that the matrix $A_{m+1-k}$ is extreme in
${\mathbf{U}}^{m+1-k}.$
 Therefore, due to Proposition \ref{AalphaBbetaCalphabeta} the following matrix
\begin{eqnarray*}
A_{m+1}&=&\left(
      \begin{array}{ccc}
        A_k & \circleddash_{k\times m-k} & \circleddash_{k\times 1}\\
        \circleddash_{m-k\times k} & \circleddash_{m-k\times m-k} & \left(\frac12\right)_{m-k\times 1} \\
        \circleddash_{1\times k}& \left(\frac12\right)_{1\times m-k} & 1 \\
        \end{array}
    \right)\\
&=&\left(
      \begin{array}{cc}
        A_k & \circleddash_{k\times m+1-k} \\
        \circleddash_{m+1-k\times k} & A_{m+1-k} \\
      \end{array}
    \right)
\end{eqnarray*}
is extreme in ${\mathbf{U}}^{m+1}$, and it contains the matrix $A_m$
as a principal submatrix.

\textsc{Step II.} Now let us prove that
${\mathbf{U}}_m\hookrightarrow{\mathbf{U}}^{m+1}.$ Let
$A_m\in{\mathbf{U}}_m$ be any matrix. Then according to the
Krein-Milman theorem we have
\begin{eqnarray*}
A_m=\sum\limits_{i=1}^n\lambda_iA^{(i)}_{m},
\end{eqnarray*}
where $\lambda_i\geq0,$ $\sum\limits_{i=1}^n\lambda_i=1$ and
$A_m^{(i)}\in{\mathbf{Extr}\mathbf{U}}_m.$  Due to \texttt{I-Step},
for every $i=\overline{1,m}$ there exists a  matrix
$A^{(i)}_{m+1}\in{\mathbf{Extr}\mathbf{U}}^{m+1}$
\begin{eqnarray*}
A^{(i)}_{m+1}=\left(
      \begin{array}{cc}
        A^{(i)}_m & \left(A^{k_i,m-k_i}_{\theta,\frac12}\right)^t \\
        A^{k_i,m-k_i}_{\theta,\frac12} & 1 \\
      \end{array}
    \right),
\end{eqnarray*}
containing the matrix $A^{(i)}_{m}$ as a principal submatrix. Here
$A^{k_i,m-k_i}_{\theta,\frac12}=\left(\circleddash_{1\times
k_i},\left(\frac12\right)_{1\times m-k_i}\right),$ Then, the
following matrix
$$
A_{m+1}=\sum\limits_{i=1}^n\lambda_iA^{(i)}_{m+1},
$$
belongs to ${\mathbf{U}}^{m+1}$ and it contains the matrix $A_{m}$
as a principal submatrix. This completes the proof.
\end{pf}

\begin{cor}\label{inclusionsExtrUm}
We have the following inclusions:
\begin{itemize}
  \item [(i)] ${\mathbf{U}}^{1}\subset{\mathbf{U}}_{1}\hookrightarrow{\mathbf{U}}^{2}
\subset{\mathbf{U}}_{2}\hookrightarrow
      \cdots\hookrightarrow{\mathbf{U}}^{m-1}\subset{\mathbf{U}}_{m-1}
\hookrightarrow{\mathbf{U}}^{m}\subset{\mathbf{U}}_{m};$
  \item [(ii)]${\mathbf{Extr}\mathbf{U}}^{1}\subset{\mathbf{Extr}\mathbf{U}}_{1}
\hookrightarrow
      \cdots\hookrightarrow{\mathbf{Extr}\mathbf{U}}^{m-1}\subset{\mathbf{Extr}\mathbf{U}}_{m-1}
\hookrightarrow{\mathbf{Extr}\mathbf{U}}^{m}\subset
      {\mathbf{Extr}\mathbf{U}}_{m}.$
\end{itemize}
\end{cor}

\begin{center}
\section{Algebraic structure of the sets ${\mathbf{U}}_m$ and ${\mathbf{U}}^m$}\label{section5}
\end{center}

In this section, we are going to study an algebraic structure of the
sets ${\mathbf{U}}_m$ and ${\mathbf{U}}^m$.

Let us consider the following matrix equation
\begin{eqnarray}\label{matrixeq}
\frac{X_m+X^{t}_m}{2}=A_m,
\end{eqnarray}
where $A_m$ is a given symmetric matrix and $X_m$ is unknown matrix,
$X^{t}_m$ is  the transpose of $X_m.$

In this section, we are going to solve the following problem:
\emph{find necessary and sufficient conditions for $A_m$ in which
the matrix equation \eqref{matrixeq} has a solution in the class of
all (sub)stochastic matrices.}

We will use the following result which has been proved in \cite{GS}.

\begin{prop}\label{GS}{\rm\cite{GS}}
Let $A_m$ be a extreme matrix in ${\mathbf{U}}^m.$ If $A_m$ has no
any saturated principal submatrices of order $m-1$ then $A_m$ is a
stochastic matrix.
\end{prop}

\begin{thm}\label{stochasticcase}
Let $A_m$ be a symmetric matrix with nonnegative entries. For solvability of equation
\begin{eqnarray}\label{matrixeqstochasticcase}
\frac{X_m+X^{t}_m}{2}=A_m,
\end{eqnarray}
in the class of stochastic matrices it is necessary and sufficient
to be $A_m\in{\mathbf{U}}^m.$
\end{thm}
\begin{pf}
\emph{Necessity.} Let a stochastic matrix $X_m$ be a solution of
\eqref{matrixeqstochasticcase}. We want to show that
$A_m\in{\mathbf{U}}^m.$ Indeed, one can see that $A_m^{t}=A_m$ and
\begin{eqnarray*}
\sum\limits_{i,j\in\alpha}a_{ij}&=&\frac12\left(\sum\limits_{i,j\in\alpha}x_{ij}+\sum\limits_{i,j\in\alpha}x_{ji}\right)\\
&=&\sum\limits_{i,j\in\alpha}x_{ij}=\sum\limits_{i\in\alpha}\sum\limits_{j\in\alpha}x_{ij}\le\sum\limits_{i\in\alpha}1=|\a|,
\end{eqnarray*}
for any $\a\subset I_m.$ Moreover, if $\a=I_m$ then we have
\begin{eqnarray*}
\sum\limits_{i,j=1}^ma_{ij}&=&\frac12\left(\sum\limits_{i,j=1}^mx_{ij}+\sum\limits_{i,j=1}^mx_{ji}\right)\\
&=&\sum\limits_{i,j=1}^mx_{ij}=\sum\limits_{i=1}^m\sum\limits_{j=1}^mx_{ij}=\sum\limits_{i=1}^m1=m,
\end{eqnarray*}
which means that $A_m\in{\mathbf{U}}^m.$

\emph{Sufficiency.} Let $A_m\in{\mathbf{U}}^m.$ We must to show the
existence of a stochastic matrix $X_m$ for which
\eqref{matrixeqstochasticcase} is satisfied.

First, assume that $A_m\in{\mathbf{Extr}\mathbf{U}}^m.$ In this
case, we use induction with respect to the order of $A_m.$
Elementary calculations show that the assertion of the theorem is
true for $m=2.$ We assume that the assertion of the theorem is true
for all $m\le k-1$ and we prove it for $m=k.$

If $A_k$ has no any saturated principal submatrices of order $k-1,$
then according to Proposition \ref{GS}, $A_k$ is a stochastic
matrix. In this case, as a solution of equation
\eqref{matrixeqstochasticcase} we can take $A_k$ itself.

Let us assume $A_k$ has a saturated principal submatrix of order
$k-1.$ We denote it by $A_{\a},$ where $|\a|=k-1.$ Since
$A_k\in{\mathbf{Extr}\mathbf{U}}^k,$ due to Proposition
\ref{extrUmandextrUm}, we have $A_k\in{\mathbf{Extr}\mathbf{U}}_k.$
Since $A_{\a}$ is a saturated principal submatrix of $A_k$ and
$A_k\in{\mathbf{Extr}\mathbf{U}}_k,$ according to Theorem
\ref{localcriterion} (ii), we get
$A_{\a}\in{\mathbf{Extr}\mathbf{U}}_{|\a|}.$ From $|\a|=k-1,$ due to
the assumption of induction, for the matrix $A_\a$ there exists a
solution of equation \eqref{matrixeqstochasticcase} in the class of
stochastic matrices. We denote this solution by
$X'_\a=(x'_{ij})_{i,j\in\a}.$ Let $\{i_0\}=I_m\setminus\a.$ We know
that $A_k$ and $A_\a$ are saturated matrices, then one can get
\begin{eqnarray}\label{ai0i0plus2ai0j}
a_{i_0i_0}+2\sum\limits_{j\in\a}a_{i_0j}=1.
\end{eqnarray}
By $A_k\in{\mathbf{Extr}\mathbf{U}}_k,$ according to Theorem
\ref{localcriterion} (ii), equality \eqref{ai0i0plus2ai0j} yields
the following possible two cases

\textsc{Case I:} $a_{i_0i_0}=1$ and $a_{i_0j}=a_{ji_0}=0$ for all $j\in\a$;

\textsc{Case II:} $a_{i_0j_0}=a_{j_0i_0}=\frac12$ for some $j_0\neq i_0$ and $a_{i_0j}=a_{ji_0}=0$  for all $j\in I_m\setminus\{j_0\}.$

In \textsc{Case I}, we define a solution $X_m=(x_{ij})_{i,j=1}^m$ of
equation \eqref{matrixeqstochasticcase}, corresponding to the matrix
$A_m,$ as follows
\begin{eqnarray*}
x_{ij}=\left\{\begin{array}{l}
         x'_{ij} \quad i,j\in\a\\
         0 \quad \, \, \, i=i_0,\, \, j\in\a\\
         0 \quad \, \, \, j=i_0, \, \, i\in\a\\
         1 \quad \, \, \, i=i_0, \, \, j=i_0
       \end{array}
       \right.,
\end{eqnarray*}
where $X'_\a=(x'_{ij})_{i,j\in\a}$ is a solution of equation \eqref{matrixeqstochasticcase}, corresponding to the matrix $A_\a.$ One can easily check that $X_m$ is a stochastic matrix.

In \textsc{Case II}, let us define a solution
$X_m=(x_{ij})_{i,j=1}^m$ of equation \eqref{matrixeqstochasticcase},
corresponding to the matrix $A_m,$ as follows
\begin{eqnarray*}
x_{ij}=\left\{\begin{array}{l}
         x'_{ij} \quad i,j\in\a\\
         0 \quad \, \, \, i=i_0,\, \, j\in I_m\setminus\{j_0\}\\
         0 \quad \, \, \, j=i_0, \, \, i\in I_m\\
         1 \quad \, \, \, i=i_0, \, \, j=j_0
       \end{array}
       \right.,
\end{eqnarray*}
where $X'_\a=(x'_{ij})_{i,j\in\a}$ is a solution of equation
\eqref{matrixeqstochasticcase}, corresponding to the matrix $A_\a.$
One can easily check that $X_m$ is a stochastic matrix.

So, for extreme matrices of the set ${\mathbf{U}}^m,$ the assertion
of the theorem has been proved.

Now, we are going to prove it, for any elements of the set
${\mathbf{U}}^m.$

Let $A_m\in {\mathbf{U}}^m.$ According to Krein-Milman's theorem,
$A_m$ can be represented as the convex combination of extreme
matrices of ${\mathbf{U}}^m,$ i.e.
\begin{eqnarray}\label{matrixAm}
A_m=\sum\limits_{i=1}^s\lambda_iA^{(i)}_m,
\end{eqnarray}
where, $0\le\lambda_i\le1,$ $\sum\limits_{i=1}^s\lambda_i=1,$ and
$A^{(i)}_m\in{\mathbf{Extr}\mathbf{U}}^m$ for all
$i=\overline{1,s}.$

By $X^{(i)}_m,$ we denote solutions of equation
\eqref{matrixeqstochasticcase}, corresponding to the extreme
matrices $A^{(i)}_m$ of ${\mathbf{U}}^m$, where $i=\overline{1,s}.$

We define a matrix $X_m$ as follows
\begin{eqnarray}\label{matrixXm}
X_m=\sum\limits_{i=1}^s\lambda_iX^{(i)}_m.
\end{eqnarray}
Since every matrix $X^{(i)}_m$ is  stochastic, the matrix $X_m$
defined by \eqref{matrixXm} is a solution of equation
\eqref{matrixeqstochasticcase} in the class of stochastic matrices,
corresponding to $A_m.$
\end{pf}

By means of Theorem \ref{stochasticcase} we are going to generalize this result for substochastic matrices.

\begin{thm}\label{substochasticcase}
Let $A_m$ be a symmetric matrix with nonnegative entries. For solvability of equation
\begin{eqnarray}\label{matrixeqsubstochasticcase}
\frac{X_m+X^{t}_m}{2}=A_m,
\end{eqnarray}
in the class of substochastic matrices it is necessary and
sufficient to be $A_m\in{\mathbf{U}}_m.$
\end{thm}

\begin{pf}
\emph{Necessity.} Let a substochastic matrix $X_m$ be a solution of
equation \eqref{matrixeqsubstochasticcase}. We want to show that
$A_m\in{\mathbf{U}}_m.$ Indeed, one can see that $A_m^{t}=A_m$ and
\begin{eqnarray*}
\sum\limits_{i,j\in\alpha}a_{ij}&=&\frac12\left(\sum\limits_{i,j\in\alpha}x_{ij}+\sum\limits_{i,j\in\alpha}x_{ji}\right)\\
&=&\sum\limits_{i,j\in\alpha}x_{ij}=\sum\limits_{i\in\alpha}\sum\limits_{j\in\alpha}x_{ij}\le\sum\limits_{i\in\alpha}1=|\a|,
\end{eqnarray*}
for any $\a\subset I_m,$ this means that $A_m\in{\mathbf{U}}_m.$

\emph{Sufficiency.} Let $A_m\in{\mathbf{U}}_m.$ We must to show the
existence of a substochastic matrix $X_m$ for which
\eqref{matrixeqsubstochasticcase} is satisfied.

As above proved theorem, we shall prove the assertion, for extreme
matrices of ${\mathbf{U}}_m.$ Then, we prove it for any elements of
${\mathbf{U}}_m.$

Let $A_m\in{\mathbf{Extr}\mathbf{U}}_m.$ If $A_m$ is a saturated
matrix, then due to Proposition \ref{extrUmandextrUm},
$A_m\in{\mathbf{Extr}\mathbf{U}}^m.$ According to Theorem
\ref{stochasticcase}, there exists a solution $X_m$ of equation
\eqref{matrixeqsubstochasticcase}, corresponding to $A_m,$ in the
class of stochastic matrices. Since every stochastic matrix is
substochastic, the matrix $X_m$ is a solution of equation
\eqref{matrixeqsubstochasticcase} in the class of substochastic
matrices.

If $A_m$ is not a saturated matrix, then due to Corollary
\ref{nonsaturatedcase}, there exist two index sets $\a$ and $\a'$
with $|\a'|=K$ and $\a=I_m\setminus\a'$ such that $A_\a$ is a
saturated principal submatrix of $A_m$ containing all nonzero
entries of $A_m$, and for any $i\in\a'$ every $i^{\rm{th}}$ row and
$i^{\rm{th}}$ column of $A_m$ consists zeroes, where
$K=m-\sum\limits_{i,j=1}^ma_{ij}.$ Since $A_\a$ is a saturated there
exists a solution $X'_\a=(x'_{ij})_{i,j\in\a}$ of equation
\eqref{matrixeqsubstochasticcase} in the class of substochastic
matrices, corresponding to $A_\a.$ Now, using the matrix $X'_\a,$ we
construct a substochastic matrix $X_m=(x_{ij})_{i,j=1}^m$ as follows
\begin{eqnarray*}
x_{ij}=\left\{\begin{array}{l}
         x'_{ij} \quad i,j\in\a\\
         0 \quad \, \, \, i\in\a', \, \, j\in I_m\\
         0 \quad \, \, \, j\in\a', \, \, i\in I_m
       \end{array}
       \right.,
\end{eqnarray*}
which is a solution of equation \eqref{matrixeqsubstochasticcase} corresponding to the matrix $A_m.$

Hence, for extreme matrices of the set ${\mathbf{U}}_m,$ the
assertion of the theorem has been proved.

For any elements of the set ${\mathbf{U}}_m$ the proof can be
proceeded by the same argument as in the proof of Theorem
\ref{stochasticcase}.
\end{pf}

\section*{Acknowledgments}

The authors are grateful to Professor Yuri Safarov for his valuable
comments and remarks on improving the paper. The authors acknowledge
the MOSTI grants 01-01-08-SF0079 and CLB10-04. The second named
author (F.M.) acknowledges the Junior Associate scheme of the Abdus
Salam International Centre for Theoretical Physics, Trieste, Italy.

\end{document}